\documentclass[11pt,a4paper,twoside]{article}
\usepackage{latexsym,amsmath,amsthm,amssymb,amsfonts}
\usepackage{mathrsfs,multicol}
\usepackage{amscd}
\usepackage{cite}
\usepackage[ansinew]{inputenc}

\def\proof{{\bf Proof.}\quad}
\def\endproof{\hfill\hbox{$\sqcup$}\llap{\hbox{$\sqcap$}}\medskip}
\numberwithin{equation}{section}
\newcommand{\zdhy}{\allowdisplaybreak}

\catcode`@=11

\newskip\plaincentering \plaincentering=0pt plus 1000pt minus 1000pt
\def\@plainlign{\tabskip=0pt\everycr={}}
\def\eqalignno#1{\displ@y \tabskip\plaincentering
  \halign to\displaywidth{\hfil$\@lign\displaystyle{##}$\tabskip\z@skip
    &$\@lign\displaystyle{{}##}$\hfil\tabskip\plaincentering
    &\llap{$\@lign##$}\tabskip\z@skip\crcr
    #1\crcr}}
\def\leqalignno#1{\displ@y \tabskip\plaincentering
  \halign to\displaywidth{\hfil$\@lign\displaystyle{##}$\tabskip\z@skip
    &$\@lign\displaystyle{{}##}$\hfil\tabskip\plaincentering
    &\kern-\displaywidth\rlap{$\@lign##$}\tabskip\displaywidth\crcr
    #1\crcr}}
\def\plainLet@{\relax\iffalse{\fi\let\\=\cr\iffalse}\fi}
\def\plainvspace@{\def\vspace##1{\noalign{\vskip##1}}}

\def\intic@{\mathchoice{\hskip5\p@}{\hskip4\p@}{\hskip4\p@}{\hskip4\p@}}
\def\negintic@
 {\mathchoice{\hskip-5\p@}{\hskip-4\p@}{\hskip-4\p@}{\hskip-4\p@}}
\def\intkern@{\mathchoice{\!\!\!}{\!\!}{\!\!}{\!\!}}
\def\intdots@{\mathchoice{\cdots}{{\cdotp}\mkern1.5mu
    {\cdotp}\mkern1.5mu{\cdotp}}{{\cdotp}\mkern1mu{\cdotp}\mkern1mu
      {\cdotp}}{{\cdotp}\mkern1mu{\cdotp}\mkern1mu{\cdotp}}}
\newcount\intno@
\def\iint{\intno@=\tw@\futurelet\next\ints@}
\def\iiint{\intno@=\thr@@\futurelet\next\ints@}
\def\iiiint{\intno@=4 \futurelet\next\ints@}
\def\idotsint{\intno@=\z@\futurelet\next\ints@}
\def\ints@{\findlimits@\ints@@}
\newif\iflimtoken@
\newif\iflimits@
\def\findlimits@{\limtoken@false\limits@false\ifx\next\limits
 \limtoken@true\limits@true\else\ifx\next\nolimits\limtoken@true\limits@false
    \fi\fi}
\def\multintlimits@{\intop\ifnum\intno@=\z@\intdots@
  \else\intkern@\fi
    \ifnum\intno@>\tw@\intop\intkern@\fi
     \ifnum\intno@>\thr@@\intop\intkern@\fi\intop}
\def\multint@{\int\ifnum\intno@=\z@\intdots@\else\intkern@\fi
   \ifnum\intno@>\tw@\int\intkern@\fi
    \ifnum\intno@>\thr@@\int\intkern@\fi\int}
\def\ints@@{\iflimtoken@\def\ints@@@{\iflimits@
   \negintic@\mathop{\intic@\multintlimits@}\limits\else
    \multint@\nolimits\fi\eat@}\else
     \def\ints@@@{\multint@\nolimits}\fi\ints@@@}
\def\Sb{_\bgroup\vspace@
        \baselineskip=\fontdimen10 \scriptfont\tw@
        \advance\baselineskip by \fontdimen12 \scriptfont\tw@
        \lineskip=\thr@@\fontdimen8 \scriptfont\thr@@
        \lineskiplimit=\thr@@\fontdimen8 \scriptfont\thr@@
        \Let@\vbox\bgroup\halign\bgroup \hfil$\scriptstyle
            {##}$\hfil\cr}
\def\endSb{\crcr\egroup\egroup\egroup}
\def\Sp{^\bgroup\vspace@
        \baselineskip=\fontdimen10 \scriptfont\tw@
        \advance\baselineskip by \fontdimen12 \scriptfont\tw@
        \lineskip=\thr@@\fontdimen8 \scriptfont\thr@@
        \lineskiplimit=\thr@@\fontdimen8 \scriptfont\thr@@
        \Let@\vbox\bgroup\halign\bgroup \hfil$\scriptstyle
            {##}$\hfil\cr}
\def\endSp{\crcr\egroup\egroup\egroup}
\def\Let@{\relax\iffalse{\fi\let\\=\cr\iffalse}\fi}
\def\vspace@{\def\vspace##1{\noalign{\vskip##1 }}}
\def\aligned{\,\vcenter\bgroup\plainvspace@\plainLet@\openup\jot\m@th\ialign
  \bgroup \strut\hfil$\displaystyle{##}$&$\displaystyle{{}##}$\hfil\crcr}
\def\endaligned{\crcr\egroup\egroup}
\def\matrix{\,\vcenter\bgroup\plainLet@\plainvspace@
    \normalbaselines
  \m@th\ialign\bgroup\hfil$##$\hfil&&\quad\hfil$##$\hfil\crcr
    \mathstrut\crcr\noalign{\kern-\baselineskip}}
\def\endmatrix{\crcr\mathstrut\crcr\noalign{\kern-\baselineskip}\egroup
                \egroup\,}
\newtoks\hashtoks@
\hashtoks@={#}
\def\format{\crcr\egroup\iffalse{\fi\ifnum`}=0 \fi\format@}
\def\format@#1\\{\def\preamble@{#1}%
  \def\c{\hfil$\the\hashtoks@$\hfil}%
  \def\r{\hfil$\the\hashtoks@$}%
  \def\l{$\the\hashtoks@$\hfil}%
  \setbox\z@=\hbox{\xdef\Preamble@{\preamble@}}\ifnum`{=0 \fi\iffalse}\fi
   \ialign\bgroup\span\Preamble@\crcr}

\def\cases{\left\{\,\vcenter\bgroup\plainvspace@
     \normalbaselines\openup\jot\m@th
      \plainLet@\ialign\bgroup$\displaystyle{##}$\hfil&
      \quad$\displaystyle{{}##}$\hfil\crcr
      \mathstrut\crcr\noalign{\kern-\baselineskip}}

\newif\iftagsleft@
\tagsleft@true
\def\TagsOnRight{\global\tagsleft@false}
\def\tag#1$${\iftagsleft@\leqno\else\eqno\fi
 \hbox{\def\pagebreak{\global\postdisplaypenalty-\@M}%
 \def\nopagebreak{\global\postdisplaypenalty\@M}\rm(#1\unskip)}%
  $$\postdisplaypenalty\z@\ignorespaces}
\interdisplaylinepenalty=\@M
\def\allowdisplaybreak{\noalign{\allowbreak}}
\def\plainallowdisplaybreak@{\def\allowdisplaybreak{\noalign{\allowbreak}}}
\def\plaindisplaybreak@{\def\displaybreak{\noalign{\break}}}
\def\align#1\endalign{\def\tag{&}\plainvspace@\plainallowdisplaybreak@
\plaindisplaybreak@
  \iftagsleft@\plainlalign@#1\endalign\else
   \plainralign@#1\endalign\fi}
\def\plainralign@#1\endalign{\displ@y\plainLet@\tabskip\plaincentering
\halign to\displaywidth
     {\hfil$\displaystyle{##}$\tabskip=\z@&$\displaystyle{{}##}$\hfil
       \tabskip=\plaincentering&\llap{\hbox{\rm(##\unskip)}}\tabskip\z@\crcr
             #1\crcr}}
\def\plainlalign@
 #1\endalign{\displ@y\plainLet@\tabskip\plaincentering\halign to \displaywidth
   {\hfil$\displaystyle{##}$\tabskip=\z@&$\displaystyle{{}##}$\hfil
   \tabskip=\plaincentering&\kern-\displaywidth
        \rlap{\hbox{\rm(##\unskip)}}\tabskip=\displaywidth\crcr
               #1\crcr}}

\def\re@#1{\par\hangindent\parindent\indent\llap{#1\enspace}\ignorespaces}
\def\qfootnote#1{\edef\@sf{\spacefactor\the\spacefactor}{}#1\@sf
      \insert\footins{\let\egroup=}\footnotesize 
      \interlinepenalty100 \let\par=\endgraf
        \leftskip=0pt \rightskip=0pt
        \splittopskip=10pt plus 1pt minus 1pt \floatingpenalty=20000
   \smallskip\re@{#1}\bgroup\strut\aftergroup{\strut\egroup}\let\next}
\catcode`\@=\active
\TagsOnRight
\begin{document}
\title{\bf Gradient estimates and entropy formulae of porous medium
and fast diffusion equations for the Witten Laplacian}
\author{Guangyue Huang\thanks{Research supported by NSFC (No. 11001076, 11171091).},\
\ Haizhong Li\thanks{Research supported by NSFC (No. 10971110).}}
\date{}
\maketitle
\begin{quotation}
\noindent{\bf Abstract.} We consider gradient estimates to positive
solutions of porous medium equations and fast diffusion equations:
$$u_t=\Delta_\phi(u^p)$$
associated with the Witten Laplacian on Riemannian manifolds. Under
the assumption that the $m$-dimensional Bakry-Emery Ricci curvature
is bounded from below, we obtain gradient estimates which generalize
the results in \cite{pengni09} and \cite{Huangli}. Moreover,
inspired by X. -D. Li's work in \cite{lixiangdong11} we also study the
entropy formulae introduced in \cite{pengni09} for porous medium
equations and fast diffusion equations associated with the Witten
Laplacian. We prove monotonicity theorems for such entropy formulae
on compact Riemannian manifolds with non-negative $m$-dimensional
Bakry-Emery
Ricci curvature.\\

 \noindent{{\bf Keywords.} porous medium equation,
fast diffusion equation,
 entropy formulae, Witten Laplacian
} \\
{{\bf Mathematics Subject Classification.}  Primary 35B45, Secondary
35K55}
\end{quotation}

\section{Introduction}

Let $(M^n,g)$ be an $n$-dimensional complete Riemannian manifold.
Li and Yau \cite{{liyau86}} studied positive solutions of the heat
equation
\begin{equation}\label{Int1}
u_t=\Delta u
\end{equation} and obtained the following gradient estimates:

\noindent{\bf Theorem A(Li-Yau\cite{liyau86}).} {\it Let $(M^n,g)$
be a complete Riemannian manifold with ${\rm Ric}(B_p(2R))\\
\geq-K$, $K\geq0$. Suppose that $u$ is a positive solution of
\eqref{Int1} on $B_p(2R)\times [0,T]$. Then on $B_p(R)$,
\begin{equation}\label{Int2}
\frac{|\nabla u|^2}{u^2} - \alpha\frac{u_t}{u}
\leq\frac{C(n)\alpha^2}{R^2}\left(\frac{\alpha^2}{\alpha-1}+\sqrt{K}
R \right)+\frac{n\alpha^2K}{2(\alpha-1)}+\frac{n\alpha^2}{2t},
\end{equation} where $\alpha>1$ is a constant. Moreover, when $R\rightarrow \infty$,
\eqref{Int2} yields the following estimate on complete noncompact
Riemannian manifold $(M^n,g)$:
\begin{equation}\label{Int3}
\frac{|\nabla u|^2}{u^2}-\alpha\frac{u_t}{u} \leq
\frac{n\alpha^2K}{2(\alpha-1)}+\frac{n\alpha^2}{2t}.\end{equation} }

Recently, J. F. Li and X. J. Xu \cite{lixu11} obtained new Li-Yau type gradient
estimates for positive solutions of the heat equation \eqref{Int1}
on Riemannian manifolds. For the related research and some improvements on
Li-Yau type gradient estimates of the equation \eqref{Int1},  see
\cite{Yau94,Yau95,Bakry99,hamilton93,li05,davies89} and the
references therein. The equation
\begin{equation}\label{Int4}
u_t=\Delta(u^p)
\end{equation} with $p>1$ is called the porous medium
equation, which is a nonlinear version of the classical heat
equation. For various values of $p>1$, it has arisen in different
applications to model diffusive phenomena (see
\cite{vazquez07,Aronson79,pengni09} and the references therein). The
equation \eqref{Int4} with $p\in(0, 1)$ is called the fast diffusion
equation, which appears in plasma physics and in geometric flows.
However, there are marked differences between the porous medium
equations and the fast diffusion equation, see
\cite{Vazquez06,Daskalopoulos07}. For gradient estimates of
\eqref{Int4}, see \cite{Huangli,Aronson79,vazquez07,xuxiangjin}.

In \cite{pengni09}, Lu, Ni, V\'{a}zquez and Villani studied gradient
estimates of \eqref{Int4} and proved the following results (see
Theorem 3.3 in \cite{pengni09}):

\noindent{\bf Theorem B(P. Lu, L. Ni, J. V\'{a}zquez,
C.Villani\cite{pengni09}).} {\it Let $(M^n,g)$ be a complete
Riemannian manifold with ${\rm Ric}(B_p(2R))\geq-K$, $K\geq0$.
Suppose that $u$ is a positive solution to \eqref{Int4} with $p>1$.
Let $v=\frac{p}{p-1}u^{p-1}$ and $M=(p-1)\max_{B_p(2R)\times
[0,T]}v$. Then for any $\alpha>1$, on the ball $B_p(R)$, we have
\begin{equation}\label{maysec5}\aligned
 \frac{|\nabla v|^2}{v} -
\alpha\frac{v_t}{v}
\leq&\frac{C(n)Ma\alpha^2}{R^2}\left(\frac{\alpha^2}{\alpha-1}\frac{ap^2}{p-1}+(1+\sqrt{K}R)
\right)\\
&+\frac{\alpha^2}{\alpha-1}aMK+\frac{a\alpha^2}{t},
\endaligned\end{equation}
where $a=\frac{n(p-1)}{n(p-1)+2}$. Moreover, when $R\rightarrow
\infty$, \eqref{maysec5} yields the following estimate on complete
noncompact Riemannian manifold
$(M^n,g)$:\begin{equation}\label{maysec5-1}\aligned
 \frac{|\nabla v|^2}{v} -
\alpha\frac{v_t}{v}\leq&\frac{\alpha^2}{\alpha-1}aMK+
\frac{a\alpha^2}{t}.
\endaligned\end{equation}}

Now, we rewrite the inequality \eqref{maysec5-1} as
\begin{equation}\label{remark11}\aligned
|\nabla v|^2  - \alpha v_t\leq&\frac{\alpha^2}{\alpha-1}aMKv+
\frac{a\alpha^2v}{t}.
\endaligned\end{equation}
Since $(p-1)v=pu^{p-1}$, we have $(p-1)v\rightarrow 1$ as
$p\rightarrow 1$. Hence, $M\rightarrow 1$,
$$\aligned
|\nabla v|^2&\rightarrow \frac{|\nabla
u|^2}{u^2},\\
 v_t&\rightarrow \frac{u_t}{u},\\
 av &\rightarrow \frac{n}{2}, \endaligned$$ as $p\rightarrow 1$.
As a result, \eqref{remark11} becomes the inequality \eqref{Int3} in
Theorem A of Li-Yau. Therefore, for complete noncompact Riemannian
manifold $(M^n,g)$, the estimate \eqref{maysec5-1} in Theorem B of
Lu, Ni, V\'{a}zquez and Villani reduces to the estimate \eqref{Int3}
in Theorem A of Li-Yau when $p\rightarrow 1$.

Let $\phi\in C^2(M^n)$. The Witten Laplacian associated with $\phi$
is defined by
$$\Delta_\phi=\Delta-\nabla\phi \cdot\nabla$$ which is symmetric with respect to
the $L^2(M^n)$ inner product under the weighted measure
$$d\mu=e^{-\phi}dv,$$ that is,
$$\int\limits_{M^n}u\Delta_\phi v\,d\mu=-\int\limits_{M^n}\nabla u\nabla v\,d\mu=\int\limits_{M^n}v\Delta_\phi u\,d\mu, \ \ \ \forall \ u,v\in
C^\infty_0(M^n).$$ The $m$-dimensional
 Bakry-Emery Ricci curvature associated with the Witten Laplacian is given by
$${\rm Ric}_{\phi}^m={\rm Ric}+\nabla^2\phi-\frac{1}{m-n}d\phi\otimes
d\phi,
$$ where $m>n$ and $m=n$ if and only if $\phi$ is a constant. Define
$${\rm Ric}_{\phi}={\rm Ric}+\nabla^2\phi.$$ Then ${\rm Ric}_{\phi}$ can be seen as
the $\infty$-dimensional Bakry-Emery Ricci curvature. In this paper,
we study the following equation associated with the Witten
Laplacian:
\begin{equation}\label{Int5}
u_t=\Delta_\phi(u^p)
\end{equation} with $p>0$ and $p\neq 1$.
For $p>1$ and $p\in(0,1)$, we derive estimates of Lu,
Ni, V\'{a}zquez and Villani and Davies's type estimate. Moreover,
for $p>1$, we obtain Hamilton's type estimate and estimates of J. F. Li
and X. J. Xu. In particular, our results generalize the ones in
\cite{Huangli}.

First we consider gradient estimates of \eqref{Int5} under the
assumption that the $m$-dimensional Bakry-Emery Ricci curvature is
bounded from blew, and obtain the following results:

\noindent{\bf Theorem 1.1.} {\it Let $(M^n,g)$ be a complete
Riemannian manifold with ${\rm Ric}_{\phi}^m(B_p(2R))\geq-K$,
$K\geq0$. Suppose that $u$ is a positive solution to the porous
medium equation \eqref{Int5} with $p>1$. Let
$v=\frac{p}{p-1}u^{p-1}$ and $M=(p-1)\max_{B_p(2R)\times [0,T]}v$.
Then for any $\alpha>1$, on the ball $B_p(R)$, we have
\begin{equation}\label{Int6}\aligned
 \frac{|\nabla v|^2}{v}-
\alpha\frac{v_t}{v} \leq&\tilde{a}\alpha^2M\frac{C(m)}{R^2}
\Bigg\{\frac{\alpha^2}{\alpha-1}\frac{\tilde{a}p^2}{p-1}
+\left(1+\sqrt{K}R\coth (\sqrt{K}
R)\right)\Bigg\}\\
&+\frac{\alpha^2}{(\alpha-1)}\tilde{a}MK+\frac{\tilde{a}\alpha^2
}{t},
\endaligned\end{equation}
where $\tilde{a}=\frac{m(p-1)}{m(p-1)+2}$. Moreover, when
$R\rightarrow \infty$, \eqref{Int6} yields the following estimate on
complete noncompact Riemannian manifold
$(M^n,g)$:\begin{equation}\label{Int7}\aligned
 \frac{|\nabla v|^2}{v} -
\alpha\frac{v_t}{v}\leq&\frac{\alpha^2}{\alpha-1}\tilde{a}MK+
\frac{\tilde{a}\alpha^2}{t}.
\endaligned\end{equation}}

\noindent{\bf Theorem 1.2.} {\it Let $(M^n,g)$ be a complete
Riemannian manifold with ${\rm Ric}_{\phi}^m(B_p(2R))\geq-K$,
$K\geq0$. Suppose that $u$ is a positive solution to the fast
diffusion equation \eqref{Int5} with $p\in (1-\frac{2}{m}, 1)$. Let
$v=\frac{p}{p-1}u^{p-1}$ and $M=(1-p)\max_{B_p(2R)\times
[0,T]}(-v)$. Then for any $0<\alpha<1$, on the ball $B_p(R)$, we
have
\begin{equation}\label{Int8}\aligned
-\frac{|\nabla v|^2}{v}+\alpha\frac{v_t}{v}
\leq&\frac{(-\tilde{a})\alpha^2M}{A(\varepsilon_1,\varepsilon_2)}
\frac{C(m)}{R^2}\Bigg\{\frac{(-\tilde{a})\alpha^2p^2}{2\varepsilon_2(1-\tilde{a})(1-\alpha)(1-p)}
+\left(1+\sqrt{K}R\coth (\sqrt{K}
R)\right)\Bigg\}\\
&+\frac{(-\tilde{a})\alpha^2
MK}{\sqrt{\varepsilon_1(1-\alpha)(1-\alpha-\tilde{a})A(\varepsilon_1,\varepsilon_2)}}
+\frac{(-\tilde{a})\alpha^2}{A(\varepsilon_1,\varepsilon_2)t},
\endaligned\end{equation}
where $\tilde{a}=\frac{m(p-1)}{m(p-1)+2}$ and positive constants
$\varepsilon_1, \varepsilon_2\in (0,1)$ satisfying
$$A(\varepsilon_1,\varepsilon_2):=[1-\tilde{a}(1-\alpha)]
-\frac{(1+\varepsilon_2)^2(1-\tilde{a})^2(1-\alpha)}{(1-\varepsilon_1)(1-\alpha-\tilde{a})}
>0.$$  When $R\rightarrow
\infty$ and $\alpha\rightarrow 1$, \eqref{Int8} yields the following
estimate on complete noncompact Riemannian manifold $(M^n,g)$ with
${\rm Ric}_{\phi}^m\geq0$:
\begin{equation}\label{Int9}\aligned -\frac{|\nabla
v|^2}{v}+\frac{v_t}{v}\leq&-\frac{\tilde{a}}{t}.
\endaligned\end{equation}}

\noindent{\bf Remark 1.1.} Clearly, our estimate \eqref{Int7}
reduces to \eqref{maysec5-1} of Lu, Ni, V\'{a}zquez and Villani (see \cite{pengni09}) by
letting $m=n$. Moreover, for $p\in(0,1)$, Theorem 4.1 in
\cite{pengni09} of Lu, Ni, V\'{a}zquez and Villani can be obtained
from our Theorem 1.2 by taking $m=n$.

\noindent{\bf Theorem 1.3.} {\it Let $(M^n,g)$ be a complete
Riemannian manifold with ${\rm Ric}_{\phi}^m(B_p(2R))\geq-K$,
$K\geq0$. Suppose that $u$ is a positive solution to the porous
medium equation \eqref{Int5} with $p>1$. Let
$v=\frac{p}{p-1}u^{p-1}$ and $M=(p-1)\max_{B_p(2R)\times [0,T]}v$.
Then for any $\alpha>1$, on the ball $B_p(R)$, we have
\begin{equation}\label{thmhuang1}\aligned
 \frac{|\nabla v|^2}{v} -
\alpha\frac{v_t}{v}\leq&\tilde{a}\alpha^2\left\{\frac{\tilde{a}^{\frac{1}{2}}\alpha
p M^{\frac{1}{2}}
}{(p-1)^{\frac{1}{2}}(\alpha-1)^{\frac{1}{2}}}\frac{C(m)}{R}+\left[
\frac{1}{t}+\frac{MK}{2(\alpha -1)} \right. \right. \\
&\left.\left. +M\frac{C(m)}{R^2}\left(1+\sqrt{K}R\coth (\sqrt{K}
R)\right) \right]^{\frac{1}{2}}\right\}^2,
\endaligned\end{equation}
where $\tilde{a}=\frac{m(p-1)}{m(p-1)+2}$. Moreover, when
$R\rightarrow \infty$, \eqref{thmhuang1} yields the following
estimate on complete noncompact Riemannian manifold:
\begin{equation}\label{thmhuang2}\aligned
 \frac{|\nabla v|^2}{v} -
\alpha\frac{v_t}{v}\leq&\frac{\alpha^2}{2(\alpha-1)}\tilde{a}MK+
\frac{\tilde{a}\alpha^2}{t}.
\endaligned\end{equation}

}

\noindent{\bf Theorem 1.4.} {\it Let $(M^n,g)$ be a complete
Riemannian manifold with ${\rm Ric}_{\phi}^m(B_p(2R))\geq-K$,
$K\geq0$. Suppose that $u$ is a positive solution to the fast
diffusion equation \eqref{Int5} with $p\in (1-\frac{2}{m}, 1)$. Let
$v=\frac{p}{p-1}u^{p-1}$ and $M=(1-p)\max_{B_p(2R)\times
[0,T]}(-v)$. Then for any $0<\alpha<1$, on the ball $B_p(R)$, we
have
\begin{equation}\label{thmdavies21}\aligned
-\frac{|\nabla v|^2}{v}+\alpha\frac{v_t}{v}
\leq&\Bigg\{C(\tilde{a},\alpha)
\frac{p}{(1-p)^{\frac{1}{2}}}M^{\frac{1}{2}}\frac{C}{R}
+\Big[\Big(\frac{\alpha^2}{2(1-\alpha)}+2(1-\tilde{a})\Big)MK
+\frac{1-\alpha-\tilde{a}}{t}\\
&+(1-p)(1-\alpha-\tilde{a})M\frac{C(m)}{R^2}\Big(1+\sqrt{K}R\coth(\sqrt{K}
R)\Big)\Big]^{\frac{1}{2}}\Bigg\}^2,
\endaligned\end{equation}
where $\tilde{a}=\frac{m(p-1)}{m(p-1)+2}$. When $R\rightarrow
\infty$, \eqref{thmdavies21} yields the following estimate on
complete noncompact Riemannian manifold $(M^n,g)$:
\begin{equation}\label{thmdavies22}\aligned -\frac{|\nabla
v|^2}{v}+\alpha\frac{v_t}{v}\leq&\Big(\frac{\alpha^2}{2(1-\alpha)}+2(1-\tilde{a})\Big)MK
+\frac{1-\alpha-\tilde{a}}{t}.
\endaligned\end{equation}

}

\noindent{\bf Remark 1.2.} Our Theorem 1.3 reduces to Theorem 1.1 of
\cite{Huangli} by letting $m=n$ and the estimate \eqref{thmhuang2} improves \eqref{Int7} on complete noncompact Riemannian manifolds. For complete noncompact Riemannian
manifolds with $p\in(0,1)$, Lu, Ni, V\'{a}zquez and Villani
\cite{pengni09} proved (see Corollary 4.2 in \cite{pengni09}) the
following results: If ${\rm Ric}\geq 0$, then
\begin{equation}\label{Remark21}-\frac{|\nabla
v|^2}{v}+\frac{v_t}{v}\leq-\frac{a}{t};\end{equation} If ${\rm
Ric}\geq -K$ and $0<\alpha<1$, then  for any $\varepsilon>0$
satisfying $C(a,\alpha,\varepsilon):=1+(-a)(1-\alpha)
-\frac{(1-\alpha)(1-a)^2}{(1-\alpha)-a-(1-\alpha)\varepsilon^2}>0$,
\begin{equation}\label{Remark22}-\frac{|\nabla
v|^2}{v}+\alpha\frac{v_t}{v}\leq\frac{(-a)\alpha^2}{C(a,\alpha,\varepsilon)}
\Big(\frac{1}{t}+\frac{\sqrt{C(a,\alpha,\varepsilon)}}{(1-\alpha)\varepsilon}
MK\Big).\end{equation} Obviously, our estimate \eqref{thmdavies22}
reduces to \eqref{Remark21} of Lu, Ni, V\'{a}zquez and Villani when
$m=n$ and $\alpha\rightarrow 1$. Moreover,
\eqref{thmdavies22} is independent of $\varepsilon$.

\noindent{\bf Theorem 1.5.} {\it Let $(M^n,g)$ be a complete
Riemannian manifold with ${\rm Ric}_{\phi}^m(B_p(2R))\geq-K$,
$K\geq0$. Suppose that $u$ is a positive solution to the porous
medium equation \eqref{Int5} with $p>1$. Let
$v=\frac{p}{p-1}u^{p-1}$ and $M=(p-1)\max_{B_p(2R)\times [0,T]}v$.
Then for any $\alpha>1$, on the ball $B_p(R)$, we have
\begin{equation}\label{thhamilton1}
 \frac{|\nabla v|^2}{v}-\alpha(t)\frac{v_t}{v}
 \leq\tilde{a}\alpha^2(t)M\frac{C(m)}{R^2}\left(\frac{p^2\tilde{a}\alpha^2(t)}{2(p-1)(\alpha(t)-1)}
+3+\sqrt{K}R\coth(\sqrt{K}R)\right)+\frac{\tilde{a}\alpha^2(t)}{t},
\end{equation}
where $\tilde{a}=\frac{m(p-1)}{m(p-1)+2}$ and $\alpha(t)=e^{2MKt}$.
Moreover, when $R\rightarrow \infty$, \eqref{thhamilton1} yields the
following estimate on complete noncompact Riemannian manifold:
\begin{equation}\label{thhamilton2}\aligned
 \frac{|\nabla v|^2}{v} -
\alpha(t)\frac{v_t}{v}\leq&\frac{\tilde{a}\alpha^2(t)}{t}.
\endaligned\end{equation}
}

\noindent{\bf Remark 1.3.} Our Theorem 1.5 becomes Theorem 1.2 in
\cite{Huangli} as long as we let $m=n$.

\noindent{\bf Theorem 1.6.} {\it Let $(M^n,g)$ be a complete
Riemannian manifold with ${\rm Ric}_{\phi}^m(B_p(2R))\geq-K$,
$K\geq0$. Suppose that $u$ is a positive solution to the porous
medium equation \eqref{Int5} with $p>1$. Let
$v=\frac{p}{p-1}u^{p-1}$ and $M=(p-1)\max_{B_p(2R)\times [0,T]}v$.
Then on the ball $B_p(R)$, we have
\begin{equation}\label{maysec6}
 \frac{|\nabla v|^2}{v} -
\alpha(t)\frac{v_t}{v}-\varphi(t)
\leq\tilde{a}M\frac{C(m)}{R^2}
\Big\{1+\sqrt{K}R\coth(\sqrt{K}R)+\frac{\tilde{a}p^2}{(p-1)\tanh(MKt)}\Big\},
\end{equation}
where $\tilde{a}=\frac{m(p-1)}{m(p-1)+2}$, $\alpha(t)$ and
$\varphi(t)$ are given by
\begin{alignat}{1}\label{maysec7}
\varphi(t)
=&\tilde{a}MK\{\coth(MKt)+1\},\\
\alpha(t)=&1+\frac{\cosh(MKt)\sinh(MKt)-MKt}{\sinh^2(MKt)}.
\end{alignat}
Moreover, when $R\rightarrow \infty$, \eqref{maysec6} yields the
following estimate on complete noncompact Riemannian manifold:
\begin{equation}\label{maysec8}
 \frac{|\nabla v|^2}{v} -
\alpha(t)\frac{v_t}{v}-\varphi(t) \leq0.
\end{equation}
}

\noindent{\bf Theorem 1.7.} {\it Let $(M^n,g)$ be a complete
Riemannian manifold with ${\rm Ric}_{\phi}^m(B_p(2R))\geq-K$,
$K\geq0$. Suppose that $u$ is a positive solution to the porous
medium equation \eqref{Int5} with $p>1$. Let
$v=\frac{p}{p-1}u^{p-1}$ and $M=(p-1)\max_{B_p(2R)\times [0,T]}v$.
Then on the ball $B_p(R)$, we have
\begin{equation}\label{maysec10}
 \frac{|\nabla v|^2}{v} -
\alpha(t)\frac{v_t}{v}-\varphi(t)
\leq\tilde{a}\alpha^2(t)M\frac{C(m)}{R^2}
\Big\{1+\sqrt{K}R\coth(\sqrt{K}R)+\frac{\tilde{a}p^2\alpha^2(t)}{(p-1)\tanh(MKt)}\Big\},
\end{equation}
where $\tilde{a}=\frac{m(p-1)}{m(p-1)+2}$, $\alpha(t)$ and
$\varphi(t)$ are given by
\begin{equation}\label{maysec11}\aligned \varphi(t)
=&\frac{\tilde{a}}{t}+\tilde{a}MK+\frac{\tilde{a}}{3}(MK)^2t,\\
\alpha(t)=&1+\frac{2}{3}MKt.
\endaligned\end{equation}
Moreover, when $R\rightarrow \infty$, \eqref{maysec6} yields the
following estimate on complete noncompact Riemannian manifold:
\begin{equation}\label{maysec12}
 \frac{|\nabla v|^2}{v} -
\alpha(t)\frac{v_t}{v}-\varphi(t) \leq0.
\end{equation}
}

\noindent{\bf Remark 1.4.}  Our Theorems 1.6 and 1.7 reduce to
Theorems 1.3 and 1.4 in \cite{Huangli} by taking $m=n$,
respectively. Moreover, when $t$ is small enough,
$\alpha(t),\varphi(t)$ defined by \eqref{maysec7} and
\eqref{maysec11} both satisfy $\alpha(t)\rightarrow 1$ and
$\varphi(t)\leq2\tilde{a}MK+\frac{\tilde{a}}{t}$. Hence, \eqref{maysec8} and \eqref{maysec12} show
\begin{equation}\label{remark4}
 \frac{|\nabla v|^2}{v} -
\alpha(t)\frac{v_t}{v} \leq 2\tilde{a}MK+\frac{\tilde{a}}{t}.
\end{equation} Clearly, for $t$ small enough, \eqref{remark4} is better than
\eqref{Int7}. Therefore, \eqref{maysec8} and \eqref{maysec12}
improve \eqref{Int7} on complete noncompact Riemannian manifolds in this sense.

Denote by $R$ the scalar curvature of the metric $g$. In
\cite{Perelman02}, Perelman introduced the $\mathcal{W}$-entropy
functional as follows:
\begin{equation}\label{Int10}
\mathcal{W}(g,f,\tau)=\int\limits_{M^n}[\tau(R+|\nabla
f|^2)+f-n]\frac{e^{-f}}{(4\pi\tau)^{\frac{n}{2}}}\, dv,
\end{equation}
where $\tau$ is a positive scale parameter and $f\in C^\infty(M^n)$
satisfies
$$\int\limits_{M^n}\frac{e^{-f}}{(4\pi\tau)^{\frac{n}{2}}}\, dv=1.$$
By \cite{Perelman02}, we know that the $\mathcal{W}$-entropy is
monotone increasing under the Ricci flow, and its critical points
are given by gradient shrinking solitons. In
\cite{nilei04,2nilei04}, Ni considered the $\mathcal{W}$-entropy for
the linear heat equation
\begin{equation}\label{Int11}
u_\tau=\Delta u
\end{equation}
on complete Riemannian manifolds. More precisely, for the
$\mathcal{W}$-entropy associated with \eqref{Int11}:
\begin{equation}\label{Int12}
\mathcal{W}(g,f,\tau)=\int\limits_{M^n}[\tau|\nabla
f|^2+f-n]\frac{e^{-f}}{(4\pi\tau)^{\frac{n}{2}}}\, dv,
\end{equation} where $u=\frac{e^{-f}}{(4\pi\tau)^{\frac{n}{2}}}$ is a
positive solution to \eqref{Int11} and $\int_{M^n}\,u\, dv=1$, Ni
\cite{nilei04} proved
\begin{equation}\label{Int13}
\frac{d}{d\tau}\mathcal{W}(g,f,\tau)=-2\int\limits_{M^n}\tau\Big(\Big|\nabla^2
f-\frac{g}{2\tau}\Big|^2+{\rm Ric}(\nabla f,\nabla f)\Big)u\, dv.
\end{equation}
In particular, if the Ricci curvature is non-negative, then
$\mathcal{W}$-entropy defined by \eqref{Int13} is monotone
non-increasing on complete Riemannian manifolds. For the research of
the monotonicity of $\mathcal{W}$-entropy to other geometric heat
flows on Riemannian manifolds, see
\cite{Kotschwar09,Ecker07,nilei04,2nilei04,pengni09}. In
\cite{lixiangdong11}, Li studied the $\mathcal{W}_m$-entropy
associated with the Witten Laplacian to the linear heat equation
\begin{equation}\label{Int14}
u_\tau=\Delta_\phi u
\end{equation} on complete Riemannian manifolds satisfying the
$\mu$-bounded geometry condition. More precisely, for the
$\mathcal{W}_m$-entropy associated with \eqref{Int14}:
\begin{equation}\label{Int15}
\mathcal{W}_m(g,f,\tau)=\int\limits_{M^n}[\tau|\nabla
f|^2+f-m]\frac{e^{-f}}{(4\pi\tau)^{\frac{m}{2}}}\, d\mu,
\end{equation} where $u=\frac{e^{-f}}{(4\pi\tau)^{\frac{m}{2}}}$ is a positive
solution to \eqref{Int14}, Li \cite{lixiangdong11} proved that if
there exist two constants $m>n$ and $K\geq0$ such that ${\rm
Ric}_{\phi}^m\geq-K$, then
\begin{equation}\label{Int16}\aligned
\frac{d}{d\tau}\mathcal{W}_m(g,f,\tau)=&-2\int\limits_{M^n}\tau\Big(\Big|\nabla^2f
-\frac{g}{2\tau}\Big|^2+{\rm Ric}_{\phi}^m(\nabla f,\nabla
f)\Big)u\,d\mu\\
&-\frac{2}{m-n}\int\limits_{M^n}\tau\Big(\nabla \phi\nabla
f+\frac{m-n}{2\tau}\Big)^2u\,d\mu.
\endaligned\end{equation} In particular, if the
${\rm Ric}_{\phi}^m\geq0$, then $\mathcal{W}_m(g,f,\tau)$ is
non-increasing along the heat equation \eqref{Int14}. For the study
to the Witten Laplacian associated with the $m$-dimensional
Bakry-Emery Ricci curvature on complete Riemannian manifolds, see
\cite{Wei09,Wang04,Wang97,Qi98,Qi97,Ni02,li05,Fang09,Bakry05,Bakry94,Bakry85}.
 Let $u$ be a positive solution to
\eqref{Int4}, and let $v=\frac{p}{p-1}u^{p-1}$. In \cite{pengni09},
Lu,  Ni, V\'{a}zquez and Villani introduced the following:
$$\mathcal{N}_p(g,u,t)=-t^{a}\int\limits_{M^n}uv\, dv$$ and
\begin{equation}\label{Int17}
\mathcal{W}_p(g,u,t)=\frac{d}{dt}[t\mathcal{N}_p(g,u,t)]=t^{a+1}\int\limits_{M^n}\Big(p\frac{|\nabla
v|^2}{v}-\frac{a+1}{t}\Big)uv\, dv,
\end{equation} where $a=\frac{n(p-1)}{n(p-1)+2}$. They proved that if $M^n$ is
compact, then
\begin{equation}\label{Int18}\aligned
\frac{d}{dt}\mathcal{W}_p(g,u,t)=&-2(p-1)t^{a+1}
\int\limits_{M^n}\Big(\Big|\nabla^2v+\frac{g}{[n(p-1)+2]t}\Big|^2+{\rm
Ric}(\nabla v,\nabla v)\Big)uv\, dv\\
&-2t^{a+1}\int\limits_{M^n}\Big((p-1)\Delta v+\frac{a}{t}\Big)^2uv\,
dv.
\endaligned\end{equation}
In particular, if the Ricci curvature is non-negative, then the
entropy defined in \eqref{Int17} is monotone non-increasing on
compact Riemannian manifolds when $p>1$. For $p<1$, using the
Cauchy-Schwarz inequality, they proved from \eqref{Int18} that
\begin{equation}\label{Int19}\aligned
\frac{d}{dt}\mathcal{W}_p(g,u,t)\leq&-2t^{a+1}
\int\limits_{M^n}\Big[\frac{n(p-1)+1}{n(p-1)}\Big((p-1)\Delta
v+\frac{a}{t}\Big)^2\\
&+(p-1){\rm Ric}(\nabla v,\nabla v)\Big]uv\, dv.
\endaligned\end{equation} Clearly, if the Ricci curvature is non-negative
and $p\in(1-\frac{1}{n}, 1)$, then \eqref{Int19} shows that
$\frac{d}{dt}\mathcal{W}_p(g,u,t)\leq0$ and the entropy defined in
\eqref{Int17} is monotone non-increasing on compact Riemannian
manifolds.

Inspired by \cite{lixiangdong11}, in this paper we also study the
$\mathcal{W}_{p,m}$-entropy associated with the Witten Laplacian to
the equation \eqref{Int5} on compact Riemannian manifolds with $p>0$
and $p\neq 1$. First we define
\begin{equation}\label{addInt20}\mathcal{N}_{p,m}(g,u,t)
=-t^{\tilde{a}}\int\limits_{M^n}uv\, d\mu\end{equation} and the
$\mathcal{W}_{p,m}$-entropy is defined by
\begin{equation}\label{Int21}
\mathcal{W}_{p,m}(g,u,t)=\frac{d}{dt}[t\mathcal{N}_{p,m}(g,u,t)],
\end{equation} where $\tilde{a}=\frac{m(p-1)}{m(p-1)+2}$.
Under the $m$-dimensional Bakry-Emery Ricci curvature is bounded
from below, we prove the following:

\noindent{\bf Theorem 1.8.} {\it Let $(M^n,g)$ be a compact
Riemannian manifold. If $u$ is a positive solution to the porous
medium equation \eqref{Int5} with $p>1$, then
\begin{equation}\label{Int23}\frac{d}{dt}\mathcal{N}_{p,m}(g,u,t)
=-t^{\tilde{a}}\int\limits_{M^n}\left((p-1)\Delta_\phi
v+\frac{\tilde{a}}{t}\right)uv\,d\mu,\end{equation} where
$v=\frac{p}{p-1}u^{p-1}$ and $\tilde{a}=\frac{m(p-1)}{m(p-1)+2}$. In particular, if ${\rm
Ric}_{\phi}^m\geq0$, then
$\frac{d}{dt}\mathcal{N}_{p,m}(g,u,t)\leq0$ and
$\mathcal{N}_{p,m}(g,u,t)$ is monotone non-increasing in $t$.
Moreover,
\begin{equation}\label{addInt23}
\mathcal{W}_{p,m}(g,u,t)=t^{\tilde{a}+1}\int\limits_{M^n}\Big(p\frac{|\nabla
v|^2}{v}-\frac{\tilde{a}+1}{t}\Big)uv\,d\mu
\end{equation}
and
\begin{equation}\label{Int24}\aligned\frac{d}{dt}\mathcal{W}_{p,m}(g,u,t)
=&-2(p-1)t^{\tilde{a}+1}\int\limits_{M^n}
\Bigg\{\Big|\nabla^2v+\frac{g}{[m(p-1)+2]t}\Big|^2\\
&+\frac{1}{m-n} \Big|\nabla\phi\nabla
v-\frac{m-n}{[m(p-1)+2]t}\Big|^2+{\rm Ric}_{\phi}^m(\nabla v,\nabla
v)\Bigg\}uv\,d\mu\\
&-2t^{\tilde{a}+1}\int\limits_{M^n}\Big|(p-1)\Delta_\phi
v+\frac{\tilde{a}}{t}\Big|^2uv\,d\mu.
\endaligned\end{equation}
In particular, if ${\rm Ric}_{\phi}^m\geq0$, then
$\frac{d}{dt}\mathcal{W}_{p,m}(g,u,t)\leq0$ and
$\mathcal{W}_{p,m}(g,u,t)$ is monotone non-increasing in $t$.}

\noindent{\bf Theorem 1.9.} {\it If $u$ is a positive solution to
the fast diffusion equation \eqref{Int5} with $p\in (0, 1)$, then
\begin{equation}\label{Int25}\frac{d}{dt}\mathcal{N}_{p,m}(g,u,t)
=-t^{\tilde{a}}\int\limits_{M^n}\left((p-1)\Delta_\phi
v+\frac{\tilde{a}}{t}\right)uv\,d\mu,\end{equation} where
$v=\frac{p}{p-1}u^{p-1}$ and $\tilde{a}=\frac{m(p-1)}{m(p-1)+2}$. In particular, if ${\rm
Ric}_{\phi}^m\geq0$ and $p\in (1-\frac{2}{m}, 1)$, then
$\frac{d}{dt}\mathcal{N}_{p,m}(g,u,t)\leq0$ and
$\mathcal{N}_{p,m}(g,u,t)$ is monotone non-increasing in $t$.
Moreover,
\begin{equation}\label{addInt25}
\mathcal{W}_{p,m}(g,u,t)=t^{\tilde{a}+1}\int\limits_{M^n}\Big(p\frac{|\nabla
v|^2}{v}-\frac{\tilde{a}+1}{t}\Big)uv\,d\mu
\end{equation}
and for any positive constant $\varepsilon\geq m-n$ and
$1-\frac{1}{n+\varepsilon}\leq p\leq1-\frac{m-n}{m\varepsilon}$,
\begin{equation}\label{Int26}\aligned\frac{d}{dt}\mathcal{W}_{p,m}(g,u,t)
\leq&2t^{\tilde{a}+1}\int\limits_{M^n}\Bigg\{(1-p){\rm
Ric}_{\phi}^m(\nabla
v,\nabla v)\\
&+\Big(\frac{1-n(1-p)}{n(1-p)}-\frac{\varepsilon}{n}\Big)\Big|(p-1)\Delta_\phi v
+\frac{\tilde{a}}{t}\Big|^2\\
&+\Big(\frac{m(1-p)}{n(m-n)}-\frac{1}{n\varepsilon}\Big)
\Big|\nabla\phi\nabla
v-\frac{m-n}{[m(p-1)+2]t}\Big|^2\Bigg\}uv\,d\mu.
\endaligned\end{equation}
In particular, if ${\rm Ric}_{\phi}^m\geq0$, then
$\frac{d}{dt}\mathcal{W}_{p,m}(g,u,t)\leq0$ and
$\mathcal{W}_{p,m}(g,u,t)$ is monotone non-increasing in $t$.}

\noindent{\bf Remark 1.5.} In particular, if $m=n$, then we have
that $\phi$ is a constant. Then \eqref{Int24} becomes (5.6) of Lu,
Ni, V\'{a}zquez and Villani in \cite{pengni09}. By letting $m=n$ and  $\varepsilon\rightarrow
0$, \eqref{Int26} becomes \eqref{Int19}, which is Corollary 5.10 in
\cite{pengni09}.

\noindent{\bf Acknowledgements.} The authors would like to thank
Professor Xiang-Dong Li for his valuable comments. The authors would like to thank Mr. Zhijie Huang for helpful discussions.

\section{Proofs of Theorem 1.1 and 1.2}

Let $v=\frac{p}{p-1}u^{p-1}$. By virtue of the equation
\eqref{Int5}, we have $v_t=(p-1)v\Delta_\phi v+|\nabla v|^2$ which
is equivalent to
\begin{equation}\label{Proof1}
\frac{v_t}{v}=(p-1)\Delta_\phi v+\frac{|\nabla v|^2}{v}.
\end{equation}

\noindent{\bf Lemma 2.1.} {\it As in \cite{pengni09}, we introduce
the following differential operator
$$\mathcal{L}=\partial_t-(p-1)v\Delta_\phi.$$ Let $F=\frac{|\nabla
v|^2}{v}-\alpha\frac{v_t}{v}-\varphi$, where $\alpha=\alpha(t)$ and
$\varphi=\varphi(t)$ are functions depending on t.

(1) If $p>1$, then
\begin{equation}\label{Proof2}\aligned
\mathcal{L}(F)\leq&-\frac{1}{\tilde{a}}[(p-1)\Delta_\phi v]^2-2(p-1)
{\rm Ric}_{\phi}^m(\nabla v, \nabla v)+2p\nabla v\nabla F\\
&+(1-\alpha)\left(\frac{v_t}{v}\right)^2-\alpha'\frac{v_t}{v}-\varphi';
\endaligned\end{equation}

(2) If $p\in (0, 1)$, then
\begin{equation}\label{Proof3}\aligned
\mathcal{L}(F)\geq&-\frac{1}{\tilde{a}}[(p-1)\Delta_\phi v]^2-2(p-1)
{\rm Ric}_{\phi}^m(\nabla v, \nabla v)+2p\nabla v\nabla F\\
&+(1-\alpha)\left(\frac{v_t}{v}\right)^2-\alpha'\frac{v_t}{v}-\varphi',
\endaligned\end{equation} where $\tilde{a}=\frac{m(p-1)}{m(p-1)+2}$.

}

\proof We only give the proof to the case that $p>1$. The proof to
$p<1$ is similar, so we omit it here.

By a direct calculation, we have
\begin{equation}\label{Proof4}\mathcal{L}\Big(\frac{f}{g}\Big)=\frac{1}{g}\mathcal{L}(f)
-\frac{f}{g^2}\mathcal{L}(g)+2(p-1)v\nabla\Big(\frac{f}{g}\Big)
\nabla\log g,\ \ \ \ \  \forall \ f,g\in C^\infty(M).\end{equation}
Using \eqref{Proof1} we obtain
\begin{equation}\label{Proof5}\mathcal{L}(v_t)=(p-1)v_t\Delta_\phi v+2\nabla v \nabla v_t.\end{equation}
It is well known that for the $m$-dimensional Bakry-Emery Ricci
curvature, we have the following Bochner formula (for the elementary
proof, see\cite{Ledoux00,li05}): $$\aligned
\frac{1}{2}\Delta_\phi(|\nabla
w|^2)=&|\nabla^2w|^2+\nabla w\nabla \Delta_\phi w+{\rm Ric}_{\phi}(\nabla w, \nabla w)\\
\geq&\frac{1}{n}|\Delta w|^2+\nabla w\nabla \Delta_\phi w+{\rm Ric}_{\phi}(\nabla w, \nabla w)\\
\geq&\frac{1}{m}|\Delta_\phi w|^2+\nabla w\nabla \Delta_\phi w+{\rm
Ric}_{\phi}^m(\nabla w, \nabla w).
\endaligned$$ It follows from $p>1$ that
\begin{equation}\label{Proof6}\aligned
\mathcal{L}(|\nabla v|^2)\leq&2\nabla v \nabla
v_t-2(p-1)v\Big(\frac{1}{m}|\Delta_\phi v|^2+\nabla
v\nabla\Delta_\phi v
+{\rm Ric}_{\phi}^m(\nabla v,\nabla v)\Big)\\
=&2\nabla v\nabla [(p-1)v\Delta_\phi v+|\nabla
v|^2]-2(p-1)v\Big(\frac{1}{m}|\Delta_\phi v|^2\\
&+\nabla v \nabla \Delta_\phi v+{\rm Ric}_{\phi}^m(\nabla v,\nabla v)\Big)\\
=&2(p-1)|\nabla v|^2\Delta_\phi v+2\nabla
v\nabla(|\nabla v|^2)-\frac{2(p-1)}{m}v(\Delta_\phi v)^2\\
&-2(p-1)v{\rm Ric}_{\phi}^m(\nabla v,\nabla v).
\endaligned\end{equation}
Applying \eqref{Proof5} and \eqref{Proof6} into \eqref{Proof4}
yields
$$\mathcal{L}\left(\frac{v_t}{v}\right)=(p-1)\frac{v_t}{v}\Delta_\phi
v+\frac{2}{v}\nabla v\nabla v_t-\frac{v_t}{v}\frac{|\nabla
v|^2}{v}+2(p-1)v\nabla\left(\frac{v_t}{v}\right)\nabla(\log v),
$$
$$\aligned \mathcal{L}\left(\frac{|\nabla v|^2}{v}\right)
\leq&2(p-1)\frac{|\nabla v|^2}{v}\Delta_\phi v+\frac{2}{v}\nabla
v\nabla(|\nabla
v|^2)-\frac{2(p-1)}{m}(\Delta_\phi v)^2\\
&-2(p-1){\rm Ric}_{\phi}^m(\nabla v,\nabla v)-\frac{|\nabla
v|^4}{v^2} +2(p-1)v\nabla\left(\frac{|\nabla
v|^2}{v}\right)\nabla(\log v)\endaligned$$ and hence
\begin{equation}\aligned\label{Proof7}
\mathcal{L}(F)=&\mathcal{L}\left(\frac{|\nabla
v|^2}{v}\right)-\alpha\mathcal{L}\left(\frac{v_t}{v}\right)-\alpha'\frac{v_t}{v}-\varphi'\\
\leq&2(p-1)\frac{|\nabla v|^2}{v}\Delta_\phi v+\frac{2}{v}\nabla
v\nabla(|\nabla
v|^2)-\frac{2(p-1)}{m}(\Delta_\phi v)^2\\
&-2(p-1){\rm Ric}_{\phi}^m(\nabla v,\nabla v)-\frac{|\nabla
v|^4}{v^2} +2(p-1)v\nabla\left(\frac{|\nabla v|^2}{v}\right)\nabla(\log v)\\
&-\alpha(p-1)\frac{v_t}{v}\Delta_\phi v-\alpha\frac{2}{v}\nabla
v\nabla v_t+\alpha\frac{v_t}{v}\frac{|\nabla
v|^2}{v}-2\alpha(p-1)v\nabla\left(\frac{v_t}{v}\right)\nabla(\log
v)\\
&-\alpha'\frac{v_t}{v}-\varphi'.
\endaligned\end{equation} Noticing
$$2(p-1)v\nabla(\frac{|\nabla v|^2}{v})\nabla(\log
v)-2\alpha(p-1)v\nabla(\frac{v_t}{v})\nabla(\log v)=2(p-1) \nabla
v\nabla F,$$
$$\frac{2}{v}\nabla v\nabla(|\nabla v|^2)-\alpha\frac{2}{v}\nabla
v\nabla v_t=\frac{2}{v}\nabla v\nabla[(F+\varphi)v]
=2(F+\varphi)\frac{|\nabla v|^2}{v}+2\nabla v\nabla F,
$$ we have
\begin{equation}\label{Proof8}
\aligned 2(p-1)&v\nabla\left(\frac{|\nabla
v|^2}{v}\right)\nabla(\log
v)-2\alpha(p-1)v\nabla\left(\frac{v_t}{v}\right)\nabla(\log
v)+\frac{2}{v}\nabla v\nabla(|\nabla v|^2)-\alpha\frac{2}{v}\nabla
v\nabla v_t\\
=&2p \nabla v\nabla F+2(F+\varphi)\frac{|\nabla v|^2}{v}\\
=&2p \nabla v\nabla F+2\left(\frac{|\nabla
v|^2}{v}-\alpha\frac{v_t}{v}\right)\frac{|\nabla v|^2}{v}.
\endaligned\end{equation} On the other hand, using \eqref{Proof1}
again, we have
\begin{equation}\aligned\label{Proof9}
2(p-1)&\frac{|\nabla v|^2}{v}\Delta_\phi v-\frac{|\nabla v|^4}{v^2}
-\alpha(p-1)\frac{v_t}{v}\Delta_\phi
v+\alpha\frac{v_t}{v}\frac{|\nabla
v|^2}{v}\\
=&2\frac{|\nabla v|^2}{v}\left(\frac{v_t}{v}-\frac{|\nabla
v|^2}{v}\right)-\frac{|\nabla
v|^4}{v^2}-\alpha\frac{v_t}{v}\left(\frac{v_t}{v}-\frac{|\nabla
v|^2}{v}\right)+\alpha\frac{v_t}{v}\frac{|\nabla
v|^2}{v}\\
=&(2\alpha+2)\frac{v_t}{v}\frac{|\nabla v|^2}{v}-3\frac{|\nabla
v|^4}{v^2}-\alpha\left(\frac{v_t}{v}\right)^2.
\endaligned\end{equation}
Combining \eqref{Proof8} with \eqref{Proof9} gives
\begin{equation}\aligned\label{Proof10}
2(p-1)&v\nabla\left(\frac{|\nabla v|^2}{v}\right)\nabla(\log
v)-2\alpha(p-1)v\nabla\left(\frac{v_t}{v}\right)\nabla(\log
v)+\frac{2}{v}\nabla v\nabla(|\nabla v|^2)\\
&-\alpha\frac{2}{v}\nabla v\nabla v_t+2(p-1)\frac{|\nabla
v|^2}{v}\Delta_\phi v-\frac{|\nabla v|^4}{v^2}
-\alpha(p-1)\frac{v_t}{v}\Delta_\phi v
+\alpha\frac{v_t}{v}\frac{|\nabla v|^2}{v}\\
=&2p \nabla v\nabla F-\left(\frac{v_t}{v}- \frac{|\nabla
v|^2}{v}\right)^2+(1-\alpha)\left(\frac{v_t}{v}\right)^2\\
=&2p\nabla v\nabla F-[(p-1)\Delta_\phi
v]^2+(1-\alpha)\left(\frac{v_t}{v}\right)^2.
\endaligned\end{equation}
Putting \eqref{Proof10} into \eqref{Proof7} yields
$$\aligned
\mathcal{L}(F)\leq&-\frac{2(p-1)}{m}(\Delta_\phi v)^2
-2(p-1){\rm Ric}_{\phi}^m(\nabla v,\nabla v)+2p \nabla v\nabla F\\
&-[(p-1)\Delta_\phi
v]^2+(1-\alpha)(\frac{v_t}{v})^2-\alpha'\frac{v_t}{v}-\varphi'\\
=&-\frac{1}{\tilde{a}}[(p-1)\Delta_\phi v]^2-2(p-1){\rm Ric}_{\phi}^m(\nabla v,\nabla v)+2p \nabla v\nabla F\\
&+(1-\alpha)(\frac{v_t}{v})^2-\alpha'\frac{v_t}{v}-\varphi',
\endaligned$$ which completes the proof of (1) in Lemma 2.1.\endproof

\noindent{\bf Proof of Theorem 1.1.} Let $\xi$ be a cut-off function
such that $\xi(r)=1$ for $r\leq1$, $\xi(r)=0$ for $r\geq2$,
$0\leq\xi(r)\leq1$, and
$$0\geq\xi'(r)\geq-c_1\xi^{\frac{1}{2}}(r),$$
$$\xi''(r)\geq-c_2,$$ for positive constants $c_1$ and $c_2$.
Denote by $\rho(x)=d(x, p)$ the distance between $x$ and $p$ in
$M^n$. Let
$$\psi(x)=\xi\left(\frac{\rho(x)}{R}\right).$$ Making use of an argument
of Calabi \cite{Calabi57} (see also Cheng and Yau
\cite{Chengyau75}), we can assume without loss of generality that
the function $\psi$ is smooth in $B_p(2R)$. Then, we have
\begin{equation}\label{Proof11}\frac{|\nabla
\psi|^2}{\psi}\leq\frac{C}{R^2}.\end{equation} By the comparison
theorem with respect to the Witten Laplacian (see p. 1324,
\cite{li05})
$$\Delta_\phi\rho\geq\sqrt{(m-1)K}\coth\left(\sqrt{\frac{K}{m-1}}\ \rho\right),$$
we have
\begin{equation}\label{Proof12}\Delta_\phi\psi
=\frac{\xi'\Delta_\phi\rho}{R}+\frac{\xi''|\nabla \rho|^2}{R^2}
\geq-\frac{C(m)}{R^2}\left(1+\sqrt{K}R\coth(\sqrt{K} R)\right).
\end{equation}

Define $\tilde{F}=\frac{|\nabla v|^2}{v}-\alpha\frac{v_t}{v}$, where
$\alpha>1$ is a constant. Under the assumption that ${\rm
Ric}_{\phi}^m\geq-K$, \eqref{Proof2} shows that
\begin{equation}\label{Proof13}\aligned
\mathcal{L}(\tilde{F})\leq&-\frac{1}{\tilde{a}}[(p-1)\Delta_\phi
v]^2+2(p-1)K|\nabla v|^2
+2p\nabla v\nabla \tilde{F}\\
\leq&-\frac{1}{\tilde{a}}[(p-1)\Delta_\phi v]^2+2MK\frac{|\nabla
v|^2}{v} +2p\nabla v\nabla \tilde{F}.
\endaligned\end{equation}
Define $G=t\psi\tilde{F}$. Next we will apply maximum principle to
$G$ on $B_p(2 R)\times [0,T]$. Assume $G$ achieves its maximum at
the point $(x_0,s)\in B_p(2 R)\times [0,T]$ and assume $G(x_0,s)>0$
(otherwise the proof is trivial), which implies $s>0$. Then at the
point $(x_0,s)$, it holds that
$$\mathcal{L}(G)\geq 0, \ \ \
\nabla \tilde{F}=-\frac{\tilde{F}}{\psi} \nabla \psi$$ and by use of
\eqref{Proof13}, we have
\begin{equation}\label{Proof14}\aligned
0\leq& \mathcal{L}(G)=s\psi\mathcal{L}(\tilde{F})
-s(p-1)v\tilde{F}\Delta_\phi\psi-2s(p-1)v\nabla\tilde{F}\nabla\psi+\psi\tilde{F}\\
\zdhy
=&s\psi\mathcal{L}(\tilde{F})-(p-1)v\frac{\Delta_\phi\psi}{\psi}G+2(p-1)
v\frac{|\nabla\psi|^2}{\psi^2}G
 +\frac{G}{s}\\ \zdhy
\leq&s\psi\left(-\frac{1}{\tilde{a}}[(p-1)\Delta_\phi
v]^2+2MK\frac{|\nabla
v|^2}{v} +2p\nabla v\nabla \tilde{F}\right)\\
\zdhy &-(p-1)v\frac{\Delta_\phi
\psi}{\psi}G+2(p-1)v\frac{|\nabla\psi|^2}{\psi^2}G+\frac{G}{s}\\
\leq&-\frac{s\psi}{\tilde{a}}[(p-1)\Delta_\phi v]^2+2s\psi
MK\frac{|\nabla
v|^2}{v}+2\frac{p}{(p-1)^{\frac{1}{2}}}M^{\frac{1}{2}}G\frac{|\nabla
v|}{v^{\frac{1}{2}}}
\frac{|\nabla\psi|}{\psi}\\
\zdhy &-(p-1)v\frac{\Delta_\phi
\psi}{\psi}G+2(p-1)v\frac{|\nabla\psi|^2}{\psi^2}G+\frac{G}{s}.
\endaligned\end{equation}
Applying
$$[(p-1)\Delta_\phi v]^2=\frac{1}{\alpha^2}\tilde{F}^2
+\frac{2(\alpha-1)}{\alpha^2}\tilde{F}\frac{|\nabla
v|^2}{v}+\left(\frac{\alpha-1}{\alpha}\right)^2\frac{|\nabla
v|^4}{v^2}$$ into \eqref{Proof14}, we obtain
\begin{equation}\label{Proof15}\aligned
0\leq&-\frac{1}{\tilde{a}s\alpha^2}G^2
-\frac{2(\alpha-1)\psi}{\tilde{a}\alpha^2}G\frac{|\nabla v|^2}{v}
-\frac{s\psi^2}{\tilde{a}}\left(\frac{\alpha-1}{\alpha}\right)^2\frac{|\nabla
v|^4}{v^2}\\
&+2s\psi^2 MK\frac{|\nabla
v|^2}{v}+2\frac{p}{(p-1)^{\frac{1}{2}}}M^{\frac{1}{2}}\psi^{\frac{1}{2}}G
\frac{|\nabla v|}{v^{\frac{1}{2}}}\frac{|\nabla \psi|}{\psi^{\frac{1}{2}}}\\
&-(p-1)v (\Delta_\phi \psi) G+2(p-1)v\frac{|\nabla \psi|^2}{\psi
}G+\frac{\psi G}{s}.
\endaligned\end{equation} By virtue of the inequality
$-Ax^2+Bx\leq\frac{B^2}{4A}$, we have
$$-\frac{s\psi^2}{\tilde{a}}\left(\frac{\alpha-1}{\alpha}\right)^2\frac{|\nabla
v|^4}{v^2}+2s\psi^2 MK\frac{|\nabla
v|^2}{v}\leq\frac{\tilde{a}\alpha^2 s\psi^2M^2K^2}{(\alpha-1)^2},$$
$$-\frac{2(\alpha-1)\psi}{\tilde{a}\alpha^2}G\frac{|\nabla v|^2}{v}
+2\frac{p}{(p-1)^{\frac{1}{2}}}M^{\frac{1}{2}}\psi^{\frac{1}{2}}G
\frac{|\nabla v|}{v^{\frac{1}{2}}}\frac{|\nabla
\psi|}{\psi^{\frac{1}{2}}}\leq
\frac{\tilde{a}\alpha^2p^2M}{2(p-1)(\alpha-1)}\frac{|\nabla
\psi|^2}{\psi}G.$$ Hence, \eqref{Proof15} yields
\begin{equation}\label{Proof16}\aligned
0\leq&-\frac{1}{\tilde{a}s\alpha^2}G^2+\frac{\tilde{a}\alpha^2
s\psi^2M^2K^2}{(\alpha-1)^2}+\frac{\tilde{a}\alpha^2p^2M}{2(p-1)(\alpha-1)}\frac{|\nabla
\psi|^2}{\psi}G\\
&-(p-1)v(L\psi) G+2(p-1)v\frac{|\nabla\psi|^2}{\psi
}G+\frac{\psi G}{s}\\
\leq&-\frac{1}{\tilde{a}s\alpha^2}G^2+\Bigg\{\frac{\tilde{a}
\alpha^2p^2M}{2(p-1)(\alpha-1)}\frac{C}{R^2}
+(p-1)M\frac{C(m)}{R^2}\left(1+\sqrt{K}R\coth (\sqrt{K}
R)\right)+\frac{\psi }{s}\Bigg\}G\\
&+\frac{\tilde{a}\alpha^2 s\psi^2M^2K^2}{(\alpha-1)^2}.
\endaligned\end{equation}
Solving the quadratic inequality of $G$ in
\eqref{Proof16} yields
$$\aligned
G\leq&\frac{\tilde{a}s\alpha^2}{2}\Bigg\{\Big[\frac{\tilde{a}
\alpha^2p^2M}{2(p-1)(\alpha-1)}\frac{C}{R^2}
+M\frac{C(m)}{R^2}\left(1+\sqrt{K}R\coth (\sqrt{K}
R)\right)+\frac{\psi }{s}\Big]\\
&+\Bigg[\Big[\frac{\tilde{a}
\alpha^2p^2M}{2(p-1)(\alpha-1)}\frac{C}{R^2}
+M\frac{C(m)}{R^2}\left(1+\sqrt{K}R\coth (\sqrt{K}
R)\right)+\frac{\psi }{s}\Big]^2\\
&+\frac{4\psi^2M^2K^2}{(\alpha-1)^2} \Bigg]^{\frac{1}{2}}\Bigg\}\\
\leq&\tilde{a}s\alpha^2\Bigg\{\frac{\tilde{a}\alpha^2p^2M}{2(p-1)(\alpha-1)}\frac{C}{R^2}
+M\frac{C(m)}{R^2}\left(1+\sqrt{K}R\coth (\sqrt{K}
R)\right)+\frac{\psi }{s}+\frac{\psi M K }{(\alpha-1) } \Bigg\}.
\endaligned$$
Hence we have
\begin{equation}\label{Proof17}\aligned
G(x,T)\leq&G(x_0,s)\\
\leq&\tilde{a}T\alpha^2\frac{C(m)}{R^2}\Bigg\{\frac{\alpha^2}{(p-1)(\alpha-1)}\tilde{a}p^2M
+M\left(1+\sqrt{K}R\coth (\sqrt{K}
R)\right)\Bigg\}\\
&+\frac{\alpha^2}{(\alpha-1)}\tilde{a}TMK+\tilde{a}\alpha^2.
\endaligned\end{equation}
For all $x \in B_p(R)$, from \eqref{Proof17}, it holds that
$$\aligned F(x,T)
\leq&\tilde{a}\alpha^2M\frac{C(m)}{R^2}\Bigg\{\frac{\alpha^2}{\alpha-1}\frac{\tilde{a}p^2}{p-1}
+\left(1+\sqrt{K}R\coth (\sqrt{K}
R)\right)\Bigg\}\\
&+\frac{\alpha^2}{(\alpha-1)}\tilde{a}MK+\frac{\tilde{a}\alpha^2
}{T}.
\endaligned$$ Since $T$ is arbitrary, we complete the proof of Theorem 1.1.

\noindent{\bf Proof of Theorem 1.2.}  When $p\in (0, 1)$ we have
$v<0$ and from \eqref{Proof3}
\begin{equation}\label{Proof18}\aligned
\mathcal{L}(-\tilde{F})\leq&\frac{1}{\tilde{a}}[(p-1)\Delta_\phi
v]^2+2(p-1)
{\rm Ric}_{\phi}^m(\nabla v, \nabla v)+2p\nabla v\nabla(-\tilde{F})\\
&-(1-\alpha)\left(\frac{v_t}{v}\right)^2\\
\leq&\frac{1}{\tilde{a}}[(p-1)\Delta_\phi v]^2
+2MK\frac{|\nabla v|^2}{-v}+2p\nabla v\nabla(-\tilde{F})\\
&-(1-\alpha)\left(\frac{v_t}{v}\right)^2.
\endaligned\end{equation}
Define $G=t\psi(-\tilde{F})$. Next we will apply maximum principle
to $G$ on $B_p(2 R)\times [0,T]$. Assume $G$ achieves its maximum at
the point $(x_0,s)\in B_p(2 R)\times [0,T]$ and assume $G(x_0,s)>0$
(otherwise the proof is trivial), which implies $s>0$. Then at the
point $(x_0,s)$, it holds that
$$\mathcal{L}(G)\geq 0, \ \ \
\nabla (-\tilde{F})=-\frac{-\tilde{F}}{\psi} \nabla\psi$$ and by use
of \eqref{Proof18}, we have
\begin{equation}\label{Proof19}\aligned
0\leq& \mathcal{L}(G)
=s\psi\mathcal{L}(-\tilde{F})-(p-1)v\frac{\Delta_\phi\psi}{\psi}G+2(p-1)
v\frac{|\nabla\psi|^2}{\psi^2}G
 +\frac{G}{s}\\ \zdhy
\leq&s\psi\left(\frac{1}{\tilde{a}}[(p-1)\Delta_\phi
v]^2+2MK\frac{|\nabla v|^2}{-v}
+2p\nabla v\nabla(-\tilde{F})\right)\\
\zdhy &-(p-1)v\frac{\Delta_\phi
\psi}{\psi}G+2(p-1)v\frac{|\nabla\psi|^2}{\psi^2}G+\frac{G}{s}
-(1-\alpha)s\psi\left(\frac{v_t}{v}\right)^2\\
\leq&\frac{s\psi}{\tilde{a}}[(p-1)\Delta_\phi v]^2+2s\varphi
MK\frac{|\nabla
v|^2}{-v}+2\frac{p}{(1-p)^{\frac{1}{2}}}M^{\frac{1}{2}}G\frac{|\nabla
v|}
{(-v)^{\frac{1}{2}}}\frac{|\nabla\psi|}{\psi}\\
\zdhy &-(p-1)v\frac{\Delta_\phi
\psi}{\psi}G+2(p-1)v\frac{|\nabla\psi|^2}{\psi^2}G+\frac{G}{s}
-(1-\alpha)s\psi\left(\frac{v_t}{v}\right)^2.
\endaligned\end{equation}
Applying
$$[(p-1)\Delta_\phi v]^2=\frac{1}{\alpha^2}\tilde{F}^2
+\frac{2(\alpha-1)}{\alpha^2}\tilde{F}\frac{|\nabla
v|^2}{v}+\left(\frac{\alpha-1}{\alpha}\right)^2\frac{|\nabla
v|^4}{v^2},$$
$$\left(\frac{v_t}{v}\right)^2=
\frac{1}{\alpha^2}\left(-\tilde{F}+\frac{|\nabla
v|^2}{v}\right)^2=\frac{1}{\alpha^2}(-\tilde{F})^2
+\frac{2}{\alpha^2}(-\tilde{F})\frac{|\nabla
v|^2}{v}+\frac{1}{\alpha^2}\frac{|\nabla v|^4}{v^2}$$ into
\eqref{Proof19}, we obtain
\begin{equation}\label{Proof20}\aligned
0\leq&\frac{1}{\tilde{a}s\alpha^2}\Big\{[1-\tilde{a}(1-\alpha)]G^2
-2(1-\tilde{a})(1-\alpha)s\psi G\frac{|\nabla
v|^2}{-v}\\
&+s^2\psi^2(1-\alpha)(1-\alpha-\tilde{a})\frac{|\nabla
v|^4}{v^2}\Big\}+2s\psi^2 MK\frac{|\nabla
v|^2}{-v}\\
&+2\frac{p}{(1-p)^{\frac{1}{2}}}M^{\frac{1}{2}}\psi^{\frac{1}{2}}G
\frac{|\nabla v|}{(-v)^{\frac{1}{2}}}\frac{|\nabla
\psi|}{\psi^{\frac{1}{2}}}-(p-1)v (\Delta_\phi \psi) G\\
&+2(p-1)v\frac{|\nabla \psi|^2}{\psi }G+\frac{\psi G}{s}.
\endaligned\end{equation}
Next we take the similar method as in Theorem 4.1 of
\cite{pengni09}. Since $p\in (1-\frac{2}{m}, 1)$, we have
$\tilde{a}<0$. Thus, we have for any positive constants
$\varepsilon_1, \varepsilon_2$,
$$2s\psi^2 MK\frac{|\nabla
v|^2}{-v}\leq-\varepsilon_1\frac{s^2\psi^2}{\tilde{a}s\alpha^2}
(1-\alpha)(1-\alpha-\tilde{a})\frac{|\nabla
v|^4}{v^2}-\frac{1}{\varepsilon_1}\frac{\tilde{a}s\alpha^2(p-1)^2\psi^2
M^2K^2}{(1-\alpha)(1-\alpha-\tilde{a})},
$$
$$\aligned
2\frac{p}{(1-p)^{\frac{1}{2}}}M^{\frac{1}{2}}\psi^{\frac{1}{2}}G
\frac{|\nabla v|}{(-v)^{\frac{1}{2}}}\frac{|\nabla
\psi|}{\psi^{\frac{1}{2}}}
\leq&-\varepsilon_2\frac{2}{\tilde{a}s\alpha^2}(1-\tilde{a})(1-\alpha)s\psi
G\frac{|\nabla
v|^2}{-v}\\
&-\frac{\tilde{a}\alpha^2p^2M}{2\varepsilon_2(1-\tilde{a})(1-\alpha)(1-p)}
\frac{|\nabla\psi|^2}{\psi}G. \endaligned$$ Hence, we get from
\eqref{Proof20} that
\begin{equation}\label{Proof21}\aligned
0\leq&-\frac{1}{\tilde{a}s\alpha^2}\Big\{-[1-\tilde{a}(1-\alpha)]G^2
+2(1+\varepsilon_2)(1-\tilde{a})(1-\alpha)s\psi G\frac{|\nabla
v|^2}{-v}\\
&-(1-\varepsilon_1)s^2\psi^2(1-\alpha)(1-\alpha-\tilde{a})\frac{|\nabla
v|^4}{v^2}\Big\}-\frac{1}{\varepsilon_1}\frac{as\alpha^2\psi^2
M^2K^2}{(1-\alpha)(1-\alpha-\tilde{a})}\\
&-\frac{\tilde{a}\alpha^2p^2M}{2\varepsilon_2(1-\tilde{a})(1-\alpha)(1-p)}
\frac{|\nabla \psi|^2}{\psi}G-(p-1)v (\Delta_\phi \psi) G\\
&+2(p-1)v\frac{|\nabla \psi|^2}{\psi }G+\frac{\psi
G}{s}\\
\leq&\frac{1}{\tilde{a}s\alpha^2}\Big\{[1-\tilde{a}(1-\alpha)]
-\frac{(1+\varepsilon_2)^2(1-\tilde{a})^2
(1-\alpha)}{(1-\varepsilon_1)(1-\alpha-\tilde{a})}\Big\}G^2
-\frac{1}{\varepsilon_1}\frac{\tilde{a}s\alpha^2\psi^2
M^2K^2}{(1-\alpha)(1-\alpha-\tilde{a})}\\
&-\frac{\tilde{a}\alpha^2p^2M}{2\varepsilon_2(1-\tilde{a})(1-\alpha)(1-p)}
\frac{|\nabla \psi|^2}{\psi}G-(p-1)v (\Delta_\phi \psi) G\\
&+2(p-1)v\frac{|\nabla\psi|^2}{\psi }G+\frac{\psi G}{s}.
\endaligned\end{equation}
Taking $\varepsilon_1,\varepsilon_2$ such that
\begin{equation}\label{Proof23}[1-\tilde{a}(1-\alpha)]
-\frac{(1+\varepsilon_2)^2(1-\tilde{a})^2(1-\alpha)}{(1-\varepsilon_1)(1-\alpha-\tilde{a})}
:=A(\varepsilon_1,\varepsilon_2)>0,
\end{equation} then \eqref{Proof21} yields
\begin{equation}\label{Proof24}\aligned
0\leq&-\frac{1}{(-\tilde{a})s\alpha^2}A(\varepsilon_1,\varepsilon_2)G^2
+\Big\{\frac{(-\tilde{a})\alpha^2p^2M}{2\varepsilon_2(1-\tilde{a})(1-\alpha)(1-p)}
\frac{C}{R^2}\\
&+M\frac{C(m)}{R^2}\left(1+\sqrt{K}R\coth (\sqrt{K}
R)\right)+\frac{\psi}{s}\Big\}G+\frac{(-\tilde{a})s\alpha^2\psi^2
M^2K^2}{\varepsilon_1(1-\alpha)(1-\alpha-\tilde{a})}.
\endaligned\end{equation}
Solving the quadratic inequality of $G$ in \eqref{Proof24} yields
$$\aligned
G\leq&\frac{(-\tilde{a})s\alpha^2}{A(\varepsilon_1,\varepsilon_2)}
\Bigg\{\frac{(-\tilde{a})\alpha^2p^2M}{2\varepsilon_2(1-\tilde{a})(1-\alpha)(1-p)}
\frac{C}{R^2}+M\frac{C(m)}{R^2}\left(1+\sqrt{K}R\coth
(\sqrt{K} R)\right)+\frac{\psi}{s}\\
&+\frac{\psi M
K}{\sqrt{\varepsilon_1(1-\alpha)(1-\alpha-\tilde{a})}}
\sqrt{A(\varepsilon_1,\varepsilon_2)}\Bigg\}.
\endaligned$$
Hence we have
\begin{equation}\label{Proof25}\aligned
G(x,T)\leq&G(x_0,s)\\
\leq&\frac{(-\tilde{a})T\alpha^2M}{A(\varepsilon_1,\varepsilon_2)}
\frac{C(m)}{R^2}\Bigg\{\frac{(-\tilde{a})\alpha^2p^2}{2\varepsilon_2(1-\tilde{a})(1-\alpha)(1-p)}
+\left(1+\sqrt{K}R\coth (\sqrt{K}
R)\right)\Bigg\}\\
&+\frac{(-\tilde{a})T\alpha^2
MK}{\sqrt{\varepsilon_1(1-\alpha)(1-\alpha-\tilde{a})A(\varepsilon_1,\varepsilon_2)}}
+\frac{(-\tilde{a})\alpha^2}{A(\varepsilon_1,\varepsilon_2)}.
\endaligned\end{equation}
and for $x\in B_p(R)$,
$$\aligned -F(x,t)\leq&\frac{(-\tilde{a})\alpha^2M}{A(\varepsilon_1,\varepsilon_2)}
\frac{C(m)}{R^2}\Bigg\{\frac{(-\tilde{a})\alpha^2p^2}{2\varepsilon_2(1-\tilde{a})(1-\alpha)(1-p)}
+\left(1+\sqrt{K}R\coth (\sqrt{K}
R)\right)\Bigg\}\\
&+\frac{(-\tilde{a})\alpha^2
MK}{\sqrt{\varepsilon_1(1-\alpha)(1-\alpha-\tilde{a})A(\varepsilon_1,\varepsilon_2)}}
+\frac{(-\tilde{a})\alpha^2}{A(\varepsilon_1,\varepsilon_2)t}.
\endaligned$$
This completes the proof of Theorem 1.2.

\section{Proofs of Theorem 1.3-1.7}

Under the assumption that ${\rm Ric}_{\phi}^m\geq-K$ and $p>1$,
\eqref{Proof2} shows that
\begin{equation}\label{3Sec1}\aligned
\mathcal{L}(F)\leq&-\frac{1}{\tilde{a}}[(p-1)\Delta_\phi
v]^2+2(p-1)K|\nabla v|^2
+2p\nabla v\nabla F \\
&+(1-\alpha)\left(\frac{v_t}{v}\right)^2-\alpha'\frac{v_t}{v}-\varphi'\\
\leq&-\frac{1}{\tilde{a}}[(p-1)\Delta_\phi v]^2+2MK\frac{|\nabla
v|^2}{v} +2p\nabla v\nabla  F\\
&+(1-\alpha)\left(\frac{v_t}{v}\right)^2-\alpha'\frac{v_t}{v}-\varphi'.
\endaligned\end{equation}
Following the methods in \cite{Huangli}, we can prove that Theorem
1.3, 1.5, 1.6, 1.7 hold respectively.

Next we are in a position to prove Theorem 1.4. Define
$\overline{F}=\frac{|\nabla v|^2}{v}-\alpha\frac{v_t}{v}$, where
$0<\alpha<1$ is a constant. Then \eqref{Proof3} shows that
\begin{equation}\label{3Sec2}\aligned
\mathcal{L}(-\overline{F})\leq&\frac{1}{\tilde{a}}[(p-1)\Delta_\phi
v]^2+2MK\frac{|\nabla v|^2}{-v} +2p\nabla
v\nabla(-\overline{F})-(1-\alpha)\left(\frac{v_t}{v}\right)^2\\
=&\frac{1}{\tilde{a}\alpha^2}\Big(-\overline{F}-(1-\alpha)\frac{|\nabla
v|^2}{-v}\Big)^2+2MK\frac{|\nabla v|^2}{-v} +2p\nabla
v\nabla(-\overline{F})\\
&-\frac{1-\alpha}{\alpha^2}\Big(-\overline{F}-\frac{|\nabla
v|^2}{-v}\Big)^2.
\endaligned\end{equation}

Let $G=t\psi(-\overline{F})$. We apply maximum principle to $G$ on
$B_p(2 R)\times [0,T]$ and assume that $G$ achieves its maximum at
the point $(x_0,s)\in B_p(2 R)\times [0,T]$ with $G(x_0,s)>0$
(otherwise the proof is trivial). Then at the point $(x_0,s)$, it
holds that
$$\mathcal{L}(G)\geq 0, \ \ \
\nabla (-\overline{F})=-\frac{-\overline{F}}{\psi} \nabla\psi$$ and
by use of \eqref{3Sec2}, we have $$\aligned 0\leq& \mathcal{L}(G)
=s\psi\mathcal{L}(-\overline{F})-(p-1)v\frac{\Delta_\phi\psi}{\psi}G+2(p-1)
v\frac{|\nabla\psi|^2}{\psi^2}G
 +\frac{G}{s}\\
\leq&\frac{s\psi}{\tilde{a}\alpha^2}\Big(-\overline{F}-(1-\alpha)\frac{|\nabla
v|^2}{-v}\Big)^2+2s\varphi MK\frac{|\nabla
v|^2}{-v}+2\frac{p}{(1-p)^{\frac{1}{2}}}M^{\frac{1}{2}}G\frac{|\nabla
v|}
{(-v)^{\frac{1}{2}}}\frac{|\nabla\psi|}{\psi}\\
&-\frac{1-\alpha}{\alpha^2}s\psi\Big(-\overline{F}-\frac{|\nabla
v|^2}{-v}\Big)^2-(p-1)v\frac{\Delta_\phi
\psi}{\psi}G+2(p-1)v\frac{|\nabla\psi|^2}{\psi^2}G+\frac{G}{s} .
\endaligned$$
Let $\frac{|\nabla v |^2}{-v}=\mu(-\overline{F})$ at the point
$(x_0,s)$. Then we have $\mu\geq0$ and
\begin{equation}\label{3Sec3}\aligned
0\leq&\frac{1}{\tilde{a}\alpha^2s\psi}[1-(1-\alpha)\mu]^2G^2+2\mu
MKG+\frac{2\mu^{\frac{1}{2}}}{s^{\frac{1}{2}}\psi^{\frac{1}{2}}}
\frac{p}{(1-p)^{\frac{1}{2}}}M^{\frac{1}{2}}G^{\frac{3}{2}}
\frac{|\nabla\psi|}{\psi}\\
&-\frac{1-\alpha}{\alpha^2}\frac{1}{s\psi}(1-\mu)^2G^2-(p-1)v\frac{\Delta_\phi
\psi}{\psi}G+2(p-1)v\frac{|\nabla\psi|^2}{\psi^2}G+\frac{G}{s} .
\endaligned\end{equation}
Multiplying the both sides of \eqref{3Sec3} by $\frac{s\psi}{G}$
yields
\begin{equation}\label{3Sec4}\aligned
0\leq&\frac{1}{\tilde{a}\alpha^2}[1-(1-\alpha)\mu]^2G+2\mu
MKs\psi+2\mu^{\frac{1}{2}}s^{\frac{1}{2}}
\frac{p}{(1-p)^{\frac{1}{2}}}M^{\frac{1}{2}}
\frac{|\nabla\psi|}{\psi^{\frac{1}{2}}}G^{\frac{1}{2}}\\
&-\frac{1-\alpha}{\alpha^2}(1-\mu)^2G-(p-1)sv\Delta_\phi
\psi+2(p-1)sv\frac{|\nabla\psi|^2}{\psi}+\psi\\
=&-\tilde{A}G+2\tilde{B}G^{\frac{1}{2}}+\tilde{C},
\endaligned\end{equation} where $$\tilde{A}=
\frac{1}{-\tilde{a}\alpha^2}[1-(1-\alpha)\mu]^2+\frac{1-\alpha}{\alpha^2}(1-\mu)^2,$$
$$\tilde{B}=\mu^{\frac{1}{2}}s^{\frac{1}{2}}
\frac{p}{(1-p)^{\frac{1}{2}}}M^{\frac{1}{2}}
\frac{|\nabla\psi|}{\psi^{\frac{1}{2}}},$$
$$\tilde{C}=2\mu
MKs\psi+(1-p)s(-v)\Big(-\Delta_\phi
\psi+2\frac{|\nabla\psi|^2}{\psi}\Big)+\psi.$$

It is easy to see that
\begin{equation}\label{3Sec5}\aligned
\frac{1}{\tilde{A}}=&\frac{(-\tilde{a})\alpha^2}{[1-(1-\alpha)\mu]^2
+(-\tilde{a})(1-\alpha)(1-\mu)^2}\\
=&\frac{(-\tilde{a})\alpha^2}{1+(-\tilde{a})(1-\alpha)
-2(1-\alpha)(1-\tilde{a})\mu+(1-\alpha)(1-\alpha-\tilde{a})\mu^2}\\
\leq&1-\alpha-\tilde{a},
\endaligned\end{equation}
\begin{equation}\label{3Sec6}\aligned
\frac{2\mu}{\tilde{A}}=&\frac{2(-\tilde{a})\alpha^2\mu}{1+(-\tilde{a})(1-\alpha)
-2(1-\alpha)(1-\tilde{a})\mu+(1-\alpha)(1-\alpha-\tilde{a})\mu^2}\\
\leq&\frac{(-\tilde{a})\alpha^2}{\sqrt{[1+(-\tilde{a})(1-\alpha)](1-\alpha)(1-\alpha-\tilde{a})}
-(1-\alpha)(1-\tilde{a})}\\
=&\sqrt{[\frac{1}{1-\alpha}+(-\tilde{a})](1-\alpha-\tilde{a})}
+(1-\tilde{a})\\
\leq&\frac{\alpha^2}{2(1-\alpha)}+2(1-\tilde{a}),
\endaligned\end{equation} where the last inequality used
$\sqrt{xy}\leq\frac{1}{2}(x+y)$ and there exists a constant
$C(\tilde{a},\alpha)$ such that $\frac{\mu^{\frac{1}{2}}}{\tilde{A}}\leq
C(\tilde{a},\alpha)$. From the inequality
$\tilde{A}x^2-2\tilde{B}x\leq\tilde{C}$, we have
$x\leq\frac{2\tilde{B}}{\tilde{A}}+\sqrt{\frac{\tilde{C}}{\tilde{A}}}$.
Applying this inequality  into \eqref{3Sec4} by letting
$x=G^{\frac{1}{2}}$ gives
\begin{equation}\label{3Sec7}\aligned
G^{\frac{1}{2}}\leq&C(\tilde{a},\alpha)s^{\frac{1}{2}}
\frac{p}{(1-p)^{\frac{1}{2}}}M^{\frac{1}{2}}\frac{C}{R}
+\Big[\Big(\frac{\alpha^2}{2(1-\alpha)}+2(1-\tilde{a})\Big)MKs
+1-\alpha-\tilde{a}\\
&+(1-p)(1-\alpha-\tilde{a})Ms\frac{C(m)}{R^2}\Big(1+\sqrt{K}R\coth(\sqrt{K}
R)\Big)\Big]^{\frac{1}{2}}
\endaligned\end{equation}
Hence, for $x\in B_p(R)$, we have
\begin{equation}\label{3Sec8}\aligned
-\frac{|\nabla
v|^2}{v}+\alpha\frac{v_t}{v}\leq&\Bigg\{C(\tilde{a},\alpha)
\frac{p}{(1-p)^{\frac{1}{2}}}M^{\frac{1}{2}}\frac{C}{R}
+\Big[\Big(\frac{\alpha^2}{2(1-\alpha)}+2(1-\tilde{a})\Big)MK
+\frac{1-\alpha-\tilde{a}}{t}\\
&+(1-p)(1-\alpha-\tilde{a})M\frac{C(m)}{R^2}\Big(1+\sqrt{K}R\coth(\sqrt{K}
R)\Big)\Big]^{\frac{1}{2}}\Bigg\}^2.
\endaligned\end{equation}
We complete the proof of Theorem 1.4.

\section{Proofs of Theorem 1.8 and 1.9}

\noindent{\bf Lemma 4.1.} {\it If $M^n$ is a compact Riemannian
manifold and $u$ is a positive solution to \eqref{Int5} with $p\neq
0$, then
\begin{equation}\label{Entropy1}
\frac{d}{dt}\int\limits_{M^n} uv\,d\mu=(p-1)\int\limits_{M^n}
(\Delta_\phi v)uv\,d\mu=-p\int\limits_{M^n} |\nabla v|^2u\,d\mu.
\end{equation}
}

\proof From \eqref{Proof1}, we have
$(uv)_t=vu_t+uv_t=v\Delta_\phi(u^p)+(p-1)uv\Delta_\phi v+u|\nabla
v|^2$. It follows from $\nabla(u^p)=u\nabla v$ that
$$\int\limits_{M^n} [v\Delta_\phi(u^p)+u|\nabla
v|^2]\,d\mu=\int\limits_{M^n} [-\nabla v\nabla(u^p)+u|\nabla
v|^2]\,d\mu=0.$$ Hence
$$\aligned
\frac{d}{dt}\int\limits_{M^n} uv\,d\mu=&\int\limits_{M^n}
(uv)_t\,d\mu\\
=&\int\limits_{M^n}[v\Delta_\phi(u^p)+(p-1)uv\Delta_\phi v+u|\nabla
v|^2]\,d\mu\\
=&(p-1)\int\limits_{M^n} (\Delta_\phi v)uv\,d\mu\\
=&p\int\limits_{M^n} (\Delta_\phi v)u^p\,d\mu\\
=&-p\int\limits_{M^n} \nabla v \nabla(u^p)\,d\mu\\
=&-p\int\limits_{M^n} |\nabla v|^2u\,d\mu.
\endaligned$$
We complete the proof of Lemma 4.1.\endproof

\noindent{\bf Lemma 4.2.} {\it If $M^n$ is a compact Riemannian
manifold and $u$ is a positive solution to \eqref{Int5} with $p\neq
0$, then
\begin{equation}\label{Entropy2}
\frac{d}{dt}\int\limits_{M^n} (\Delta_\phi
v)uv\,d\mu=2\int\limits_{M^n} \Big[(p-1)(\Delta_\phi
v)^2+|\nabla^2v|^2+{\rm Ric}_\phi(\nabla v,\nabla v)\Big]uv\,d\mu.
\end{equation}
}

\proof Noticing
\begin{equation}\label{Entropy3}\frac{d}{dt}\int\limits_{M^n}
(\Delta_\phi v)uv\,d\mu=\int\limits_{M^n} [(\Delta_\phi
v)_tuv+(\Delta_\phi v)(uv)_t]\,d\mu.\end{equation} A direct
calculation gives
$$\aligned(\Delta_\phi v)_t=&\Delta_\phi[(p-1)v\Delta_\phi v+|\nabla v|^2]\\
=&(p-1)[(\Delta_\phi v)^2
+2\nabla v\nabla \Delta_\phi v+v\Delta_\phi^2v]+\Delta_\phi|\nabla v|^2\\
=&(p-1)(\Delta_\phi v)^2+2p\nabla v\nabla
\Delta_\phi v+(p-1)v\Delta_\phi^2v+2[|\nabla^2v|^2+{\rm Ric}_\phi(\nabla
v,\nabla v)].
\endaligned$$
We derive from $(p-1)\nabla(uv^2)=(2p-1)uv\nabla v$ that
$$\aligned
\int\limits_{M^n}& [2p\nabla v\nabla \Delta_\phi
v+(p-1)v\Delta_\phi^2v]uv\,d\mu\\
=&\int\limits_{M^n}2p\nabla v\nabla
(\Delta_\phi v)uv\,d\mu-\int\limits_{M^n}(p-1)\nabla(uv^2)\nabla \Delta_\phi v\,d\mu\\
=&\int\limits_{M^n}\nabla v \nabla(\Delta_\phi v)uv\,d\mu.
\endaligned$$
Hence,
\begin{equation}\label{Entropy4} \int\limits_{M^n}
(\Delta_\phi v)_tuv\,d\mu=\int\limits_{M^n} \Big\{(p-1)(\Delta_\phi
v)^2+\nabla v\nabla \Delta_\phi v +2[|\nabla^2v|^2+{\rm
Ric}_\phi(\nabla v,\nabla v)]\Big\}uv\,d\mu.
\end{equation}
On the other hand,
\begin{equation}\label{Entropy5}\aligned
\int\limits_{M^n} \Delta_\phi v(uv)_t\,d\mu=&\int\limits_{M^n}
\Delta_\phi v[v\Delta_\phi (u^p)+(p-1)uv\Delta_\phi v+u|\nabla v|^2]\,d\mu \\
=&\int\limits_{M^n}[-\nabla(v\Delta_\phi
v)\nabla(u^p)+(p-1)uv(\Delta_\phi v)^2+u|\nabla
v|^2\Delta_\phi v]\,d\mu\\
=&\int\limits_{M^n}[-\nabla(v\Delta_\phi v)u\nabla
v+(p-1)uv(\Delta_\phi v)^2+u|\nabla
v|^2\Delta_\phi v]\,d\mu\\
=&\int\limits_{M^n}[-\nabla v\nabla \Delta_\phi v +(p-1)(\Delta_\phi
v)^2]uv\,d\mu.\endaligned\end{equation} Inserting \eqref{Entropy4}
and \eqref{Entropy5} into \eqref{Entropy3} concludes the proof of
Lemma 4.2.\endproof

\noindent{\bf Proof of Theorem 1.8 and 1.9.} By Lemma 4.1, we have
$$\aligned
\frac{d}{dt}\mathcal{N}_{p,m}(g,u,t)
=&-\tilde{a}t^{\tilde{a}-1}\int\limits_{M^n}uv\,d\mu-(p-1)t^{\tilde{a}}\int\limits_{M^n}(\Delta_\phi
v)uv\,d\mu\\
=&-t^{\tilde{a}}\int\limits_{M^n}\left((p-1)\Delta_\phi
v+\frac{\tilde{a}}{t}\right)uv\,d\mu.
\endaligned$$ We obtain \eqref{Int23} and \eqref{Int25}. On the other hand,
from the definition of $\mathcal{W}_{p,m}(g,u,t)$ in \eqref{Int21},
we have
$$\aligned
\mathcal{W}_{p,m}(g,u,t)=&\frac{d}{dt}[t\mathcal{N}_{p,m}(g,u,t)]\\
=&\mathcal{N}_{p,m}(g,u,t)+t\frac{d}{dt}\mathcal{N}_{p,m}(g,u,t)\\
=&t^{\tilde{a}+1}\int\limits_{M^n}\Big(p\frac{|\nabla
v|^2}{v}-\frac{\tilde{a}+1}{t}\Big)uv\,d\mu,
\endaligned$$ where the Lemma 4.1 was used in the last equality.
Hence, we derive \eqref{addInt23} and \eqref{addInt25}.

Noticing that the estimate \eqref{Int7} also holds for compact
Riemannian manifolds. Taking $K=0$ and then letting
$\alpha\rightarrow 1$  in \eqref{Int7} yields
$$(p-1)\Delta_\phi v+\frac{\tilde{a}}{t}=\frac{v_t}{v}- \frac{|\nabla v|^2}{v}+
\frac{\tilde{a}}{t}\geq0,$$ which concludes that if ${\rm
Ric}_{\phi}^m\geq0$, then $\frac{d}{dt}\mathcal{N}_{p,m}(g,u,t)\leq0$
and $\mathcal{N}_{p,m}(g,u,t)$ is a monotone non-increasing in $t$.
When $p\in(1-\frac{2}{m}, 1)$ and ${\rm Ric}_{\phi}^m\geq0$, we also
get from \eqref{Int9} that
$$(p-1)\Delta_\phi v+\frac{\tilde{a}}{t}=\frac{v_t}{v}- \frac{|\nabla v|^2}{v}+
\frac{\tilde{a}}{t}\leq0,$$ which shows that
$\frac{d}{dt}\mathcal{N}_{p,m}(g,u,t)\leq0$ and
$\mathcal{N}_{p,m}(g,u,t)$ is also a monotone non-increasing in $t$.

Next we are in a position to prove \eqref{Int24}. From
\eqref{Int23}, we have
$$\aligned \frac{d}{dt}[t\frac{d}{dt}&\mathcal{N}_{p,m}(g,u,t)]\\
=&\frac{d}{dt}\Big[-t^{\tilde{a}+1}\int\limits_{M^n}(p-1)(\Delta_\phi
v)uv\,d\mu -\tilde{a}t^{\tilde{a}}\int\limits_{M^n}uv\,d\mu
\Big]\\
=&\frac{d}{dt}\Big[-t^{\tilde{a}+1}\int\limits_{M^n}(p-1)(\Delta_\phi
v)uv\,d\mu+\tilde{a}\mathcal{N}_{p,m}(g,u,t)
\Big]\\
=&-2t^{\tilde{a}+1}\int\limits_{M^n}\Big[(p-1)^2(\Delta_\phi
v)^2+(p-1)|\nabla^2v|^2+(p-1){\rm Ric}_\phi(\nabla v,\nabla
v)\Big]uv\,d\mu\\
&-(\tilde{a}+1)t^{\tilde{a}}\int\limits_{M^n} (p-1)(\Delta_\phi v)uv\,d\mu
-\tilde{a}t^{\tilde{a}}\int\limits_{M^n}\left((p-1)\Delta_\phi
v+\frac{\tilde{a}}{t}\right)uv\,d\mu,
\endaligned$$
where the last equality used the Lemma 4.2. Hence,
\begin{equation}\label{Entropy6}\aligned
\frac{d}{dt}&\mathcal{W}_{p,m}(g,u,t)\\
=&\frac{d}{dt}[t\frac{d}{dt}\mathcal{N}_{p,m}(g,u,t)+\mathcal{N}_{p,m}(g,u,t)]\\
=&-2t^{\tilde{a}+1}\int\limits_{M^n}\Big[(p-1)^2(\Delta_\phi
v)^2+(p-1)|\nabla^2v|^2+(p-1){\rm Ric}_\phi(\nabla v,\nabla
v)\Big]uv\,d\mu\\
&-(\tilde{a}+1)t^{\tilde{a}}\int\limits_{M^n}(p-1)(\Delta_\phi v)uv\,d\mu
-(\tilde{a}+1)t^{\tilde{a}}\int\limits_{M^n}\left((p-1)\Delta_\phi v+\frac{\tilde{a}}{t}\right)uv\,d\mu\\
=&-2t^{\tilde{a}+1}\int\limits_{M^n}\Big[(p-1)^2(\Delta_\phi
v)^2+(p-1)|\nabla^2v|^2+(p-1){\rm Ric}_\phi(\nabla v,\nabla
v)\\
&+(p-1)\frac{\tilde{a}+1}{t}\Delta_\phi v+\frac{\tilde{a}^2+\tilde{a}}{2t^2}\Big]uv\,d\mu.
\endaligned\end{equation}
Noticing $$\aligned &(p-1)^2(\Delta_\phi v)^2+(p-1)\frac{\tilde{a}+1}{t}\Delta_\phi v
+\frac{\tilde{a}^2+\tilde{a}}{2t^2}\\
=&\Big|(p-1)\Delta_\phi
v+\frac{m(p-1)}{[m(p-1)+2]t}\Big|^2+\frac{2(p-1)}{[m(p-1)+2]t}\Delta_\phi
v+\frac{(p-1)m}{[m(p-1)+2]^2t^2},
\endaligned$$
and hence \begin{equation}\label{Entropy7}\aligned &(p-1)^2(\Delta_\phi
v)^2+(p-1)\frac{\tilde{a}+1}{t}\Delta_\phi
v+\frac{\tilde{a}^2+\tilde{a}}{2t^2}+(p-1)|\nabla^2v|^2
+\frac{p-1}{m-n}(\nabla\phi\nabla v)^2\\
=&\Big|(p-1)\Delta_\phi
v+\frac{m(p-1)}{[m(p-1)+2]t}\Big|^2\\
&+(p-1)\Big|\nabla^2v+\frac{g}{[m(p-1)+2]t}\Big|^2+\frac{p-1}{m-n}
\Big|\nabla\phi\nabla v-\frac{m-n}{[m(p-1)+2]t}\Big|^2.
\endaligned\end{equation} We complete the proof of \eqref{Int24} by putting
\eqref{Entropy7} into \eqref{Entropy6}.

When $p\in (0,1)$, by the Cauchy-Schwarz inequality, we have
$$\aligned
-(p-1)&\Big|\nabla^2v+\frac{g}{[m(p-1)+2]t}\Big|^2\\
\geq&-\frac{p-1}{n}\Big|\Delta
v+\frac{n}{[m(p-1)+2]t}\Big|^2\\
=&-\frac{1}{n(p-1)}\Big|(p-1)\Delta_\phi v+\frac{\tilde{a}}{t}\Big|^2
-\frac{p-1}{n}\Big|\nabla\phi\nabla v-\frac{m-n}{[m(p-1)+2]t}\Big|^2\\
&-\frac{2}{n}\Big((p-1)\Delta_\phi
v+\frac{\tilde{a}}{t}\Big)\Big(\nabla\phi\nabla v-\frac{m-n}{[m(p-1)+2]t}\Big).
\endaligned$$
Hence,
\begin{equation}\label{Entropy8}\aligned
-(p-1)&\Big|\nabla^2v+\frac{g}{[m(p-1)+2]t}\Big|^2-\frac{p-1}{m-n}
\Big|\nabla\phi\nabla
v-\frac{m-n}{[m(p-1)+2]t}\Big|^2\\
&-\Big|(p-1)\Delta_\phi
v+\frac{\tilde{a}}{t}\Big|^2\\
\geq&\frac{1-n(1-p)}{n(1-p)}\Big|(p-1)\Delta_\phi
v+\frac{\tilde{a}}{t}\Big|^2+\frac{m(1-p)}{n(m-n)}
\Big|\nabla\phi\nabla v-\frac{m-n}{[m(p-1)+2]t}\Big|^2\\
&-\frac{2}{n}\Big((p-1)\Delta_\phi
v+\frac{\tilde{a}}{t}\Big)\Big(\nabla\phi\nabla
v-\frac{m-n}{[m(p-1)+2]t}\Big)\\
\geq&\Big(\frac{1-n(1-p)}{n(1-p)}-\frac{\varepsilon}{n}\Big)\Big|(p-1)\Delta_\phi v
+\frac{\tilde{a}}{t}\Big|^2\\
&+\Big(\frac{m(1-p)}{n(m-n)}-\frac{1}{n\varepsilon}\Big)
\Big|\nabla\phi\nabla v-\frac{m-n}{[m(p-1)+2]t}\Big|^2,
\endaligned\end{equation} where $\varepsilon\geq m-n$ is a positive
constant and satisfies $1-\frac{1}{n+\varepsilon}\leq
p\leq1-\frac{m-n}{m\varepsilon}$. Inserting
\eqref{Entropy8} into \eqref{Int24} gives
\begin{equation}\label{Entropy9}\aligned\frac{d}{dt}\mathcal{W}_{p,m}(g,u,t)
\leq&2t^{a+1}\int\limits_{M^n}\Bigg\{(1-p){\rm Ric}_{\phi}^m(\nabla
v,\nabla v)\\
&+\Big(\frac{1-n(1-p)}{n(1-p)}-\frac{\varepsilon}{n}\Big)\Big|(p-1)\Delta_\phi v
+\frac{\tilde{a}}{t}\Big|^2\\
&+\Big(\frac{m(1-p)}{n(m-n)}-\frac{1}{n\varepsilon}\Big)
\Big|\nabla\phi\nabla v-\frac{m-n}{[m(p-1)+2]t}\Big|^2\Bigg\}uv\,d\mu.
\endaligned\end{equation}
Therefore, we complete the proof of \eqref{Int26}.

\begin{flushleft}
\medskip\noindent
Guangyue Huang, \\
{Department of Mathematics, Henan Normal
University, Xinxiang 453007, P.R. China,} \\
E-mail address:~\textsf{hgy$@$henannu.edu.cn }

\medskip

Haizhong Li, \\
{Department of Mathematical Sciences,
Tsinghua University, Beijing 100084, P.R. China,}\\
E-mail address:~\textsf{hli$@$math.tsinghua.edu.cn}

\end{flushleft}

\end{document}